\documentclass[a4paper,12pt]{article} 
\usepackage{amssymb}
\usepackage{amscd}
\usepackage{amsmath} 
\usepackage{graphicx}
\usepackage{latexsym} 
\usepackage{enumerate}

\usepackage{theorem}
\theoremstyle{plain}
\theorembodyfont{\itshape}
\newtheorem{theorem}{Theorem}[section]
\newtheorem{proposition}[theorem]{Proposition}
\newtheorem{lemma}[theorem]{Lemma}

\theorembodyfont{\rmfamily}

\newtheorem{definition}[theorem]{Definition}
\newtheorem{remark}[theorem]{Remark}

\def\C{\mathbb{C}}
\def\N{\mathbb{N}}

\def\R{\mathbb{R}}
\def\Z{\mathbb{Z}}
\def\P{\mathbb{P}}
\def\a{\alpha}

\def\ag{\mathfrak{a}}
\def\bg{\mathfrak{b}}
\def\A{\varLambda}
\def\p{\mathfrak{p}}
\def\q{\mathfrak{q}}
\def\m{\mathfrak{m}}

\def\si{\sigma}

\def\Om{\Omega}
\def\om{\omega}
\def\ve{\varepsilon}
\def\vp{\varpi}

\def\wt{\widetilde} 
\def\sm{\setminus}
\def\carl{\circlearrowleft}
\newcommand{\proof}{{{\it Proof}.~~}}
\newcommand{\qed}{{\hbox{\rule[-2pt]{3pt}{6pt}}}}
\title{\bf Area-Preserving Surface Dynamics and \\
S. Saito's Fixed Point Formula\thanks{Mathematics 
Subject Classification: 37C25, 14E07, \, arXiv: {\tt 0710.0706} 
[math.DS] 3 Oct 2007}}
\author{Katsunori Iwasaki and Takato Uehara \\ \\
Graduate School of Mathematics, Kyushu University \\
6-10-1 Hakozaki, Higashi-ku, Fukuoka 812-8581 
Japan\thanks{E-mail addresses: {\tt iwasaki@math.kyushu-u.ac.jp} \ 
and \ {\tt t-uehara@math.kyushu-u.ac.jp}}} 
\begin{document}
\date{October 3, 2007}
\maketitle
\begin{abstract} 
We show that S. Saito's fixed point formula serves as a powerful 
tool for counting the number of isolated periodic points of 
an area-preserving surface map admitting periodic curves. 
His notion of periodic curves of types I and I\!I plays a central 
role in our discussion. 
We establish a Shub-Sullivan type result on the stability of local 
indices under iterations of the map, the finiteness of the number 
of periodic curves of type I\!I, and the absence of periodic 
curves of type I. 
Combined with these results, Saito's formula implies the existence 
of infinitely many isolated periodic points whose cardinality grows 
exponentially as period tends to infinity. 
\end{abstract} 
\section{Introduction} \label{sec:intro}
Counting the number of periodic points of a continuous map 
$f : X \to X$ on a compact manifold $X$ is an important issue 
in the field of dynamical systems. 
An efficient method for dealing with this problem is provided 
by several versions of Lefschetz fixed point formula. 
In general this formula aims at representing the Lefschetz number 
\begin{equation} \label{eqn:lef-number}
L(f) := \sum_{i} (-1)^i \, \mathrm{Tr} 
[\, f^{*} \,:\, H^i(X) \to H^i(X)\,] 
\end{equation}
in terms of some local data around the fixed points of $f$. 
In the simplest case where all the fixed points of $f$ are 
isolated in $X$, the classical Lefschetz formula states that 
\begin{equation} \label{eqn:lef-formula}
L(f) = \sum_{x \in X_0(f)} \nu_x(f), 
\end{equation}
where $X_0(f)$ is the set of all fixed points of $f$ and 
$\nu_x(f)$ is the local index of $f$ at $x \in X_0(f)$. 
\par 
In applying formula (\ref{eqn:lef-formula}) to counting periodic points, 
it is important to investigate the behavior of the indices 
$\nu_x(f^n)$ as $n$ tends to infinity, where $f^n := 
f \circ \cdots \circ f$ ($n$-times) stands for the $n$-th iterate of $f$. 
In this regard, Shub and Sullivan \cite{SS} prove the following 
theorem. 
\begin{theorem} \label{thm:SS}
Let $X$ be a smooth manifold and $f : X \to X$ a $C^1$-map. 
If $x \in X$ is an isolated fixed point of all iterates $f^n$, 
then the indices $\nu_x(f^n)$ are bounded as a function of $n \in \N$. 
\end{theorem}
Combining this theorem with the Lefschetz formula 
(\ref{eqn:lef-formula}), they also show that if the Lefschetz 
numbers $L(f^n)$ are unbounded then the set of periodic 
points of $f$ is infinite. 
The stability result as in Theorem \ref{thm:SS} is very useful for 
various problems in dynamical systems, so that the short 
paper \cite{SS} is often cited in the literature 
(for example in \cite{Brown,Favre,FL,Zhang}). 
\par
Unfortunately, the condition that the fixed points of a given map 
should be isolated is too restrictive, because one often encounters 
a map having a higher dimensional fixed point set. 
To cover this situation, many authors have made various 
generalizations of the classical Lefschetz formula 
(\ref{eqn:lef-formula}) and the resulting formulas have had many 
fruitful applications. 
We refer to the pioneering work of Atiyah and Bott \cite{AB} and 
some subsequent works \cite{MT,Patodi,TT} to cite only a few. 
In these generalizations, however, we have to assume that the 
induced linear map on the normal bundle to the fixed point set 
should not have eigenvalue $1$. 
This condition sometimes puts a severe restriction on the 
applicability of the formulas. 
A typical example is the case of an area-preserving map of a surface. 
Let $X$ be a surface endowed with an area form and $f : X \to X$ 
an area-preserving map. 
If $f$ admits a fixed curve $C$, then $f$ induces the identity map on 
the normal bundle to $C$ so that the above-mentioned generalizations 
do not apply to this map (see Remark \ref{rem:index} for a more 
precise discussion). 
This situation is quite unpleasant because area-preserving maps of 
surfaces constitute an important class of dynamical systems. 
Thus there should be a more suitable fixed point formula and also 
a suitable variation of the Shub-Sullivan theorem 
which fit into this class of maps. 
The aim of this article is to discuss these issues 
upon restricting our attention to algebraic surface maps over 
$\C$.  \par
At this stage we notice that S. Saito's fixed point formula 
\cite{Saito} is very appropriate for our purpose. 
It is the most desirable formula that is valid for any algebraic 
map $f : X \to X$ of a smooth projective surface $X$, where 
no restriction is put on the induced linear map on the normal 
bundle to the fixed point set. 
The success of his formula is due to the idea that the set 
$X_1(f)$ of all irreducible fixed curves of $f$ can be divided 
into two disjoint subsets, namely, into what he calls the curves 
of type I and those of type I\!I: 
\begin{equation} \label{eqn:divide}
X_1(f) = X_I(f) \amalg X_{I\!I}(f). 
\end{equation}
Then his formula expresses the Lefschetz number 
(\ref{eqn:lef-number}) in terms of suitably defined local 
indices $\nu_x(f)$ and $\nu_C(f)$ at the fixed points 
$x \in X_0(f)$ and at the fixed curves $C \in X_1(f)$, 
where different types of curves contribute to the formula in 
different ways; see formula (\ref{eqn:formula}) below. 
\par
Now we review Saito's formula more explicitly. 
Although his original formula in \cite{Saito} is stated for 
holomorphic maps, we restate it for birational maps. 
By this alteration the resulting formula gains a wider 
applicability in complex dynamics, while its proof 
remains almost the same. 
Since a birational map $f$ admits the indeterminacy set $I(f)$ 
at which $f$ is not defined, we have to restart 
with giving a proper definition of $X_0(f)$ and $X_1(f)$. 
Further we have to adapt the original definitions in 
\cite{Saito} of $X_I(f)$, $X_{I\!I}(f)$, $\nu_x(f)$ and 
$\nu_C(f)$ to the current setting. 
Leaving all these tasks in Section \ref{sec:formula}, 
we now accept that these concepts are defined properly. 
Recall also that the induced action on cohomology and 
so the Lefschetz number (\ref{eqn:lef-number}) are well 
defined for birational maps. 
Now a birational version of Saito's fixed 
point formula is stated as follows. 
\begin{theorem} \label{thm:ssaito} 
Let $X$ be a smooth projective surface and $f : X \to X$ a 
birational map different from the identity map. 
If the map $f$ satisfies the separation condition 
\begin{equation} \label{eqn:separate} 
I(f) \cap I(f^{-1}) = \emptyset, 
\end{equation}
then the Lefschetz number of $f$ is expressed as 
\begin{equation} \label{eqn:formula}
L(f) = \sum_{x \in X_0(f)} \nu_x(f) + 
\sum_{C \in X_{I}(f)} \chi_C \cdot \nu_C(f) + 
\sum_{C \in X_{I\!I}(f)} \tau_C \cdot \nu_C(f), 
\end{equation}
where $\chi_C$ is the Euler characteristic of the normalization 
of $C \in X_I(f)$ and $\tau_C$ is the self-intersection number 
of $C \in X_{I\!I}(f)$. 
\end{theorem} 
\par 
It turns out that if the induced linear map on the normal 
bundle to $C \in X_1(f)$ does not have eigenvalue $1$, then 
$C$ must be a fixed curve of type I (see Remark \ref{rem:index}). 
This clearly explains why the usual generalizations of Lefschetz 
fixed point formula are not sufficient --- they do not apply to 
the case where $f$ admits fixed curves of type I\!I. 
On the other hand, formula (\ref{eqn:formula}) is always valid, 
no matter which type of fixed curves are present. 
More remarkably, we shall see in this paper that the presence of 
curves of type I\!I plays a very crucial role in discussing a 
Shub-Sullivan type theorem and other related issues in our context. 
In this sense Saito's formula (\ref{eqn:formula}) is a very 
suitable fixed point formula for our purpose. 
\section{Main Results} \label{sec:main}
With the powerful fixed point formula (\ref{eqn:formula}) in hand, 
we proceed to stating the main results of this article. 
In what follows, unless otherwise stated explicitly, $f : X \to X$ is 
a nontrivial birational map of a smooth projective surface $X$, where 
$f$ is said to be nontrivial if it is of infinite order. 
In stating our main results, we do not assume the separation condition 
(\ref{eqn:separate}), which is only required for stating 
Theorem \ref{thm:ssaito}. \par 
Our first main theorem is a Shub-Sullivan type result in our context. 
It concerns, however, the invariance of local indices rather 
than their boundedness as in Theorem \ref{thm:SS}. 
\begin{theorem} \label{thm:main} 
For any fixed curve $C \in X_{I\!I}(f)$ of type I\!I, the indices 
$\nu_C(f^n)$ are independent of $n \in \N$. 
Similarly, for any fixed point $x \in X_0(f)$ through which at 
least one fixed curve of type I\!I passes, the indices 
$\nu_x(f^n)$ are independent of $n \in \N$. 
\end{theorem}
\begin{remark} \label{rem:main} 
If $C$ is a fixed curve of type I, then $\nu_C(f^n)$ may depend on 
$n \in \N$. 
Similarly, if a fixed point $x$ lies outside any fixed curve of 
type I\!I, then $\nu_x(f^n)$ may depend on $n \in \N$. 
See Remarks \ref{rem:ind1} and \ref{rem:ind2} for some examples 
illustrating these remarks. 
The reason why the invariance of indices is more relevant than their 
boundedness in our context is also stated there. 
\end{remark} 
\par 
Theorem \ref{thm:main} shows a special role played by curves 
of type I\!I. 
The next main theorem exhibits another role played by these curves. 
To state it we introduce a bit of terminology. 
\begin{definition} \label{def:periodic}
An irreducible curve $C$ is called a {\sl periodic} curve of $f$ 
if $C \in X_1(f^n)$ for some $n \in \N$. 
It is said to be of {\sl prime period} $n$ if $C \in X_1(f^n)$ 
but $C \not\in X_1(f^m)$ for every $m < n$. 
A periodic curve $C$ of prime period $n$ is said to be 
{\sl of type I} or {\sl of type I\!I} according as 
$C \in X_I(f^n)$ or $C \in X_{I\!I}(f^n)$, 
where the definition of $X_I(f)$ and $X_{I\!I}(f)$ is 
given later in Definition \ref{def:type2}. 
\end{definition} 
\par 
We recall some concepts from bimeromorphic surface dynamics.  
Given a bimeromorphic map $f : X \to X$ on a compact 
K\"{a}hler surface $X$, its first dynamical degree $\lambda(f)$ is 
defined by 
\begin{equation} \label{eqn:FDD}
\lambda(f) := \lim_{n \to \infty} 
|\!|(f^n)^*|_{H^{1,1}(X)} |\!|^{1/n} 
\ge 1, 
\end{equation}
where $|\!|\cdot|\!|$ is an operator norm on 
$\text{End} \, H^{1,1}(X)$. 
It is known that the limit exists, $\lambda(f)$ is independent of 
the norm $|\!| \cdot |\!|$ chosen and invariant under bimeromorphic 
conjugation (see \cite{DF}). 
The smallest possible value $\lambda(f) = 1$ corresponds to the 
case of low dynamical complexity and so we are more interested 
in the case $\lambda(f) > 1$ of higher dynamical complexity. 
A bimeromorphic map $f$ is said to be {\sl algebraically stable} 
(AS for short) if the condition $(f^n)^* = (f^*)^n : H^{1,1}(X) 
\to H^{1,1}(X)$ holds for every $n \in \N$. 
It is a standard condition under which bimeromorphic surface 
dynamics is often discussed. 
If $f$ is AS, then the first dynamical degree (\ref{eqn:FDD}) 
coincides with the spectral radius of the map $f^* : H^{1,1}(X) \to 
H^{1,1}(X)$. 
It is known that $f$ is AS if and only if 
\begin{equation} \label{eqn:AS}
f^{-m}I(f) \cap f^n I(f^{-1}) = \emptyset 
\qquad \mbox{for every} \quad 
m, n \ge 0. 
\end{equation}
\par 
Now the second main theorem is concerned with the finiteness of 
the number of periodic curves of type I\!I and also with the question: 
what happens if $f$ has `too many' periodic curves of type I\!I?
\begin{theorem} \label{thm:main2}
Let $f : X \to X$ be an AS bimeromorphic map of a compact K\"{a}hler 
surface $X$. 
\begin{enumerate}
\item If $\lambda(f) > 1$, then $f$ has at most $\rho(X) + 1$ 
irreducible periodic curves of type I\!I with mutually distinct 
prime periods, where $\rho(X)$ is the Picard number of $X$. 
\item If $\lambda(f) = 1$ and $f^n$ is not isotopic to the 
identity for any $n \in \N$, then $f$ preserves a unique rational 
or elliptic fibration $\pi : X \to S$. 
In addition, if $f$ has more than $\rho(X) + 1$ irreducible periodic 
curves of type I\!I with mutually distinct prime periods, then any 
irreducible periodic curve of type I\!I is contained in a fiber of 
the fibration $\pi$. 
\end{enumerate}
\end{theorem}
\begin{remark} \label{rem:main2} 
Two remarks are in order regarding Theorem \ref{thm:main2}. 
\begin{enumerate}
\item Assume that $f$ is nontrivial. 
Then $f$ has at most finitely many irreducible periodic curves of 
type I\!I with mutually distinct prime periods if and only if 
the number of all irreducible periodic curves of type I\!I is 
finite (see also Remark \ref{rem:P(f)}). 
\item If the Kodaira dimension of $X$ is nonnegative, then the 
bound $\rho(X)+1$ in assertion (1) of Theorem \ref{thm:main2} 
can be replaced by $\rho(X)$. 
See the proof of Theorem \ref{thm:main2} in Section \ref{sec:Fib}. 
\end{enumerate}
\end{remark}
\par 
Theorem \ref{thm:main2} asserts that the criterion $\lambda(f) > 1$ 
for a high dynamical complexity implies the finiteness of periodic 
curves of type I\!I, with an explicit bound $\rho(X) + 1$ given only 
in terms of the geometry of $X$. 
Also in the case $\lambda(f) = 1$ of low dynamical complexity, the 
presence of `too many' periodic curves of type I\!I beyond the same 
bound implies an even simpler dynamical behavior of $f$, meaning 
that any periodic curve of type I\!I is along the fibration. 
A key ingredient to establish Theorem \ref{thm:main2} is the fact that 
two periodic curves of type I\!I with distinct prime periods must be 
disjoint (see Theorem \ref{thm:prime}). 
The classification of AS bimeromorphic surface maps due to 
Diller and Favre \cite{DF} is also an important ingredient. 
\par 
The third main theorem is concerned with the intimate relationship 
between the area-preserving property of a map $f$ and the absence 
of fixed curves of type I. 
Here we interpret the area-preserving property in a wide sense 
to the effect that the area form $\omega$ preserved by $f$ 
may have poles or zeros, though of course $\omega$ must 
be nontrivial. 
\begin{theorem} \label{thm:main3}
Assume that $f : X \to X$ preserves a nontrivial meromorphic $2$-form 
$\omega$. 
If $C$ is a fixed curve along which $\omega$ has no pole 
of order $\nu_C(f)$, then $C$ must be of type I\!I. 
In particular any fixed curve is of type I\!I unless it is an 
irreducible component of the pole divisor $(\omega)_{\infty}$ of 
$\omega$. 
\end{theorem}
\begin{remark} \label{rem:pole}
If $C$ is an irreducible component of $(\omega)_{\infty}$ along which 
$\omega$ has a pole of order $\nu_C(f)$, then $C$ can be a fixed 
curve of type I. 
We refer to Remark \ref{rem:pre} for such an example. 
\end{remark}
\par 
Looking back on what we have stated, we notice that the principal 
concept underlying all the above theorems is the curves of type I\!I. 
Combined with these theorems, Saito's fixed point formula 
(\ref{eqn:formula}) is applied to the iterates of an area-preserving 
AS birational map to yield a useful formula for the number of isolated 
periodic points of it (see Theorem \ref{thm:formula}). 
This formula has an interesting implication as mentioned below. 
To state it, let $\mathrm{Per}_n^{i}(f)$ be the set of all isolated 
periodic points of $f$ with (not necessarily prime) period $n$ and 
$\# \mathrm{Per}_n^{i}(f)$ its cardinality counted with 
multiplicity (see (\ref{eqn:card-iso}) for the precise definition). 
\begin{theorem} \label{thm:main4} 
Let $f : X \to X$ be an AS birational map of a smooth projective 
surface $X$ and assume that $f$ preserves a nontrivial 
meromorphic $2$-form $\omega$ such that 
\par\smallskip 
$(*)$ \, no irreducible component of the pole divisor 
$(\omega)_{\infty}$ of $\omega$ is a periodic curve of type I. 
\par\smallskip\noindent 
If $\lambda(f) > 1$ then $f$ has at most finitely many 
irreducible periodic curves and infinitely many isolated periodic 
points. 
Moreover the number of isolated periodic points of period $n$, 
counted with multiplicity, is estimated as 
\[
|\# \mathrm{Per}_n^{i}(f) - \lambda(f)^n| \le 
\left\{ \begin{array}{ll}
O(1) \quad & (\mbox{if $X \sim$ no Abelian surface}), \\[2mm] 
4 \, \lambda(f)^{n/2} + O(1) \quad & 
(\mbox{if $X \sim$ an Abelian surface}), 
\end{array} \right. 
\]
where $O(1)$ stands for a bounded function of $n \in \N$ and 
$X \sim Y$ indicates that $X$ is birationally equivalent to $Y$. 
\end{theorem}  
In this sense, area-preserving surface maps nicely fit into 
Saito's fixed point formula (\ref{eqn:formula}) and hence the 
title of this article. 
Note that, once $f$ and $\omega$ are given concretely, the condition 
$(*)$ is verifiable in finite procedures, since $(\omega)_{\infty}$ 
contains only finitely many irreducible components. 
\par 
The plan of this article is as follows. 
After a more detailed review of Saito's fixed point formula in 
Section \ref{sec:formula}, Theorems \ref{thm:main}, 
\ref{thm:main2}, \ref{thm:main3} and \ref{thm:main4} are proved 
in Sections \ref{sec:st}, \ref{sec:Fib}, \ref{sec:APD} and 
\ref{sec:iso-per} respectively. 
Actually a refined version of Theorem \ref{thm:main4} is established 
in Section \ref{sec:iso-per} (see Theorem \ref{thm:main4dash}). 
In Section \ref{sec:example} we illustrate our main theorems 
by giving an interesting example of an area-preserving AS map on 
the minimal resolution of a singular cubic surface. 
\section{S. Saito's Fixed Point Formula} \label{sec:formula}
In this section we introduce the terminology and concepts needed to 
formulate Theorem \ref{thm:ssaito}, the birational version of 
S. Saito's fixed point formula in \cite{Saito}. 
These terminology and concepts are also needed to formulate our 
main theorems in Section \ref{sec:main}. 
\par 
We begin with some ring-theoretical preparations. 
Let $A := \C[\![z_1,z_2]\!]$ be the ring of formal power series of 
two variables over $\C$ with its maximal ideal $\m \subset A$ 
and let $\si : A \to A$ be a nontrivial continuous endomorphism 
of $A$ in the $\m$-adic topology. 
Then $\si$ is expressed as 
\begin{equation} \label{eqn:si}
\left\{
\begin{array}{rcl}
\si(z_1) &=& z_1 + g \cdot h_1, \\[1mm]
\si(z_2) &=& z_2 + g \cdot h_2, 
\end{array}
\right. 
\end{equation}
for some elements $g$, $h_1$, $h_2 \in A$, where 
$g$ is nonzero and $h_1$, $h_2$ are relatively prime. 
Consider the ideals $\ag(\si) := (g)$ and 
$\bg(\si) := (h_1, \, h_2)$ generated by $g$ and by $h_1$, $h_2$, 
respectively. 
They are independent of the choice of the coordinates $z_1$ and $z_2$. 
Since $h_1$ and $h_2$ are relatively prime, the quotient vector 
space $A/\bg(\si)$ is finite-dimensional, so that one can put 
\begin{equation} \label{eqn:delta}
\delta(\si) := \dim_{\C} \, A/\bg(\si) < \infty. 
\end{equation}
Let $\A(\si)$ be the set of all prime ideals $\p$ of height $1$ in 
$A$ that divide $\ag(\si)$. 
For $\p \in \A(\si)$, put 
\begin{equation} \label{eqn:nup}
\nu_{\p}(\si) := \max\{\, m \in \N \, | \, \ag(\si) \subset \p^m \}. 
\end{equation}
Let $\kappa[\p]$ be the normalization of the quotient ring $A/\p$, 
and $\kappa(\p)$ the quotient field of $\kappa[\p]$. 
It follows from the definition that $\kappa[\p]$ is isomorphic to 
$\C[\![t]\!]$ for some prime element $t$. 
Moreover, two modules of formal differentials are defined by 
the projective limits: 
\[
\hat{\Omega}_{A/\C}^1 := 
\lim_{\substack{\longleftarrow \\ n }} \Omega_{A_n/\C}^1, \qquad 
\hat{\Omega}_{\kappa[\p]/\C}^1 := 
\lim_{\substack{\longleftarrow \\ n }} \Omega_{\kappa[\p]_n/\C}^1, 
\]
where $A_n:=A/\m^n$ and $\kappa[\p]_n:=\C[\![t]\!]/(t)^n$. 
It is easy to see that $\hat{\Omega}_{A/\C}^1$ is a free 
$A$-module of rank two with generators $dz_1$ and $dz_2$, while 
$\hat{\Omega}_{\kappa[\p]/\C}^1$ is a free $\kappa[\p]$-module of 
rank one with a generator $dt$. 
Furthermore we put 
$\hat{\Omega}_{\kappa(\p)/\C}^1 := 
\hat{\Omega}_{\kappa[\p]/\C}^1 \otimes_{\kappa[\p]} \kappa(\p)$ 
and define a map 
\[
\tau_{\p} : \hat{\Omega}_{A/\C}^1 \to 
\hat{\Omega}_{\kappa(\p)/\C}^1 
\]
to be the homomorphism induced from the natural map 
$A \to \kappa(\p)$. 
Finally we put 
\begin{equation} \label{eqn:omega}
\vp_{\si} := 
h_2 \cdot dz_1 - h_1 \cdot dz_2 \in \hat{\Omega}_{A/\C}^1. 
\end{equation}
\begin{definition} \label{def:type} 
A prime ideal $\p \in \A(\si)$ is said to be 
{\sl of type I} relative to $\si$, if $\tau_{\p}(\vp_{\si}) \neq 0$ 
in $\hat{\Omega}_{\kappa(\p)/\C}^1$; otherwise, $\p$ is said to be 
{\sl of type I\!I}. 
Let $\A_I(\si)$ denote the set of prime ideals $\p \in \A(\si)$ of 
type I and let $\A_{I\!I}(\si)$ denote the set of prime ideals 
$\p \in \A(\si)$ of type I\!I, respectively. 
\end{definition}
\par 
In view of this definition we observe that there exists an element 
$a \in \kappa[\p]$ such that 
\begin{equation} \label{eqn:a}
\left\{
\begin{array}{rcll}
\tau_{\p}(\vp_{\si}) &=& a \cdot dt \qquad & 
(\,\mbox{if} \,\, \p \in \A_I(\si)), \\[2mm]
\vp_{\si} &=& a \cdot dp \quad \mathrm{mod} \,\,\, \p \cdot 
\hat{\Omega}_{A/\C}^1 \qquad & 
(\,\mbox{if} \,\, \p \in \A_{I\!I}(\si)), 
\end{array}
\right. 
\end{equation}
where $p \in A$ is a prime element such that $\p = (p)$. 
Identifying $\kappa[\p]$ with $\C[\![t]\!]$, we define  
\begin{eqnarray}
\mu_{\p}(\si) &:=& 
\max \{\, m \in \N \cup \{0\} \, | \, (a) \subset (t)^m \} 
\qquad (\p \in \A(\si)) \label{eqn:mup} 
\\[2mm]
\nu_{A}(\si) &:= & \delta(\si) + 
\displaystyle \sum_{\p \in \A(\si)} 
\nu_{\p}(\si) \cdot \mu_{\p}(\si). \label{eqn:nuA} 
\end{eqnarray} 
\par 
The case where $\si$ is an automorphism will be of particular 
interest later. 
\begin{lemma} \label{lem:inverse} 
If $\si$ is an automorphism, then we have:
\begin{itemize}
\item $\A(\si^{-1}) = \A(\si)$, $\A_I(\si^{-1}) = \A_I(\si)$ 
and $\A_{I\!I}(\si^{-1}) = \A_{I\!I}(\si)$$;$ 
\item $\delta(\si^{-1}) = \delta(\si)$ and 
$\nu_A(\si^{-1}) = \nu_A(\si)$$;$ 
\item $\nu_{\p}(\si^{-1}) = \nu_{\p}(\si)$ and 
$\mu_{\p}(\si^{-1}) = \mu_{\p}(\si)$ for any $\p \in \A(\si)$. 
\end{itemize} 
\end{lemma}
\proof 
As in (\ref{eqn:si}), the inverse automorphism 
$\si^{-1} : A \to A$ is expressed as 
\[
\left\{
\begin{array}{rcl}
\si^{-1}(z_1) &=& z_1 + u \cdot v_1, \\[1mm]
\si^{-1}(z_2) &=& z_2 + u \cdot v_2, 
\end{array}
\right. 
\]
for some $u$, $v_1$, $v_2 \in A$, where $u \neq 0$ and 
$v_1$, $v_2$ are relatively prime. 
Applying $\si$ to it, we have 
\[
\left\{
\begin{array}{rcl}
\si(z_1) = z_1 - \si(u) \cdot \si(v_1), \\[1mm]
\si(z_2) = z_2 - \si(u) \cdot \si(v_2), 
\end{array}
\right. 
\]
where $\si(v_1)$ and $\si(v_2)$ are relatively prime, since 
so are $v_1$ and $v_2$, and $\si$ is an automorphism. 
Comparing this with (\ref{eqn:si}) and multiplying $g$ by a 
unit and $h_i$ by its inverse, we can put $-g = \si(u)$ and 
$h_i = \si(v_i)$ for $i = 1$, $2$. 
Expanding the righthand sides into a Taylor series yields 
\begin{equation} \label{eqn:hi}
\left\{
\begin{array}{rclclcl}
-g &=& \si(u) &=& u(z_1+g \cdot h_1, \, z_2 +g \cdot h_2) 
&=& u + g \cdot b,  \\[2mm]
h_i &=& \si(v_i) &=& v_i(z_1+g \cdot h_1, \, z_2 +g \cdot h_2)
&=& v_i + g \cdot b_i \qquad (i = 1, 2),  
\end{array}
\right.
\end{equation}
with some elements $b$, $b_i \in \bg(\si)$. 
Hence $u = -g(1+b) \in \ag(\si)$ and 
$v_i = h_i - g \cdot b_i \in \bg(\si)$ for $i = 1$, $2$, 
so that $\ag(\si^{-1}) \subset \ag(\si)$ and 
$\bg(\si^{-1}) \subset \bg(\si)$. 
Replacing $\si$ with $\si^{-1}$ we also have $\ag(\si) \subset 
\ag(\si^{-1})$ and $\bg(\si) \subset \bg(\si^{-1})$. 
Thus $\ag(\si^{-1}) = \ag(\si)$ and $\bg(\si^{-1}) = \bg(\si)$.
In view of (\ref{eqn:delta}), (\ref{eqn:nup}) and the 
definition of $\A(\si)$, these equalities imply that 
$\A(\si^{-1}) = \A(\si)$, $\delta(\si^{-1}) = \delta(\si)$ and 
$\nu_{\p}(\si^{-1}) = \nu_{\p}(\si)$ for any $\p \in \A(\si)$. 
Moreover it follows from (\ref{eqn:omega}) and (\ref{eqn:hi}) 
that 
\[
\vp_{\si} = h_2 \cdot dz_1 - h_1 \cdot dz_2 
= (v_2 \cdot dz_1 - v_1 \cdot dz_2) + g 
(b_2 \cdot dz_1 - b_1 \cdot dz_2 ) \in 
\vp_{\si^{-1}} + \ag(\si) \cdot \hat{\Omega}^1_{A/\C}. 
\]
In view of (\ref{eqn:a}) and (\ref{eqn:mup}), this shows that 
$\A_I(\si^{-1}) = \A_I(\si)$, 
$\A_{I\!I}(\si^{-1}) = \A_{I\!I}(\si)$ and 
$\mu_{\p}(\si^{-1}) = \mu_{\p}(\si)$ for any $\p \in \A(\si)$. 
Finally, $\nu_A(\si^{-1}) = \nu_A(\si)$ readily follows from 
(\ref{eqn:nuA}). \qed \par\medskip
We turn our attention to birational surface maps. 
Let $f: X \to X$ be a nontrivial birational map of 
a smooth projective surface $X$. 
Since $f$ admits the indeterminacy set $I(f)$ at which $f$ 
is not defined, we must ask what should be the definition 
of $X_0(f)$, the set of fixed points of $f$. 
A natural idea is to treat the forward map $f$ and the 
backward map $f^{-1}$ symmetrically so that one can 
switch between $f$ and $f^{-1}$. 
Then even for a point $x \in I(f)$ one can declare that $x$ 
is a fixed point of $f$ provided that $x$ is away from 
$I(f^{-1})$ and is fixed by $f^{-1}$. 
The definition of $X_1(f)$ also needs some care, though it 
is just a small matter. 
\begin{definition} \label{def:fixed} 
Let $X_0^{\circ}(f)$ be the set of all points 
$x \in X \sm I(f)$ fixed by $f$ and put 
\begin{equation} \label{eqn:fixed} 
X_0(f) := X_0^{\circ}(f) \cup X_0^{\circ}(f^{-1}). 
\end{equation}
Let $X_1(f)$ be the set of all irreducible curves $C$ in $X$ 
such that $C \sm I(f)$ is fixed pointwise by $f$. 
This definition makes sense since $I(f)$ is a finite set of 
points and so $C \sm I(f)$ is a nonempty Zariski open subset 
of $C$. 
It is easy to see that the definition is symmetric: 
\begin{equation} \label{eqn:symm}
X_1(f) = X_1(f^{-1}). 
\end{equation}
\end{definition} 
\begin{remark} \label{rem:exceptional} 
Let $E(f)$ be the exceptional set of $f$. 
Any irreducible component $C$ of $E(f)$ is not an element of 
$X_1(f)$, because $C \setminus I(f)$ is contracted to a single 
point by $f$. 
\end{remark}
\par 
We now define the index $\nu_x(f)$ at a fixed point 
$x \in X_0(f)$ and the index $\nu_{C}(f)$ at a fixed curve 
$C \in X_1(f)$. 
Let $A_x$ denote the completion of the local ring of $X$ at $x$. 
If $x \in X_0^{\circ}(f)$ then the map $f$ is holomorphic 
around $x$ and hence induces a continuous endomorphism 
$f_x^* : A_x \to A_x$ with respect to the $\m_x$-adic topology 
in a natural manner, where $\m_x$ is the maximal ideal of $A_x$. 
Since $X$ is assumed to be smooth, the ring $A_x$ is isomorphic 
to the formal power series ring $\C[\![z_1,z_2]\!]$, so that upon 
putting $A = A_x$ and $\si = f_x^*$ we can come to the 
ring-theoretical situation considered above and define the 
number $\nu_{A_x}(f_x^*)$ via the formula (\ref{eqn:nuA}). 
Similarly, if $x \in X_0^{\circ}(f^{-1})$ then we can consider 
the number $\nu_{A_x}((f^{-1})_x^*)$ instead of $\nu_{A_x}(f_x^*)$. 
Moreover, if $x \in X_0^{\circ}(f) \cap X_0^{\circ}(f^{-1})$, 
then $f$ is a local biholomorphism around $x$, 
inducing an automorphism $f_x^* : A_x \to A_x$ with 
its inverse $(f_x^*)^{-1} = (f^{-1})_x^*$, so that 
Lemma \ref{lem:inverse} implies that 
\begin{equation} \label{eqn:coincidence} 
\nu_{A_x}(f_x^*) = \nu_{A_x}((f^{-1})_x^*) \qquad \mbox{at} 
\quad x \in X_0^{\circ}(f) \cap X_0^{\circ}(f^{-1}). 
\end{equation}
Next, given a fixed curve $C \in X_1(f)$, take a point 
$x$ of $C \sm I(f)$. 
Then one can speak of the continuous endomorphism 
$f_x^* : A_x \to A_x$. 
Let $C_x$ denote the germ at $x$ of the curve $C$ and let 
$\A(C_x)$ be the set of all prime ideals in $A_x$ 
determined by the irreducible components of $C_x$. 
Then any $\p \in \A(C_x)$ is a prime ideal of length $1$ that 
divides $\ag(f_x^*)$, that is, $\A(C_x) \subset \A(f_x^*)$, so 
that one can define the number $\nu_{\p}(f_x^*)$ via the formula 
(\ref{eqn:nup}) with $\si = f_x^*$. 
This definition does not depend on the choice of the point 
$x \in C \setminus I(f)$ and the ideal $\p \in \A(C_x)$ 
(see Saito \cite[page 1016]{Saito}). 
Summing up these discussions, we make the following 
definitions. 
\begin{definition} \label{def:index} 
The local index $\nu_x(f)$ at a fixed point $x \in X_0(f)$ is 
defined by 
\begin{equation} \label{eqn:nux}
\nu_x(f) := \left\{
\begin{array}{ll}
\nu_{A_x}(f_x^*) \qquad & (\,\mbox{if $x \in X_0^{\circ}(f)$}), 
\\[2mm] 
\nu_{A_x}((f^{-1})_x^*) \qquad & 
(\,\mbox{if $x \in X_0^{\circ}(f^{-1})$}), 
\end{array}\right. 
\end{equation}
which is consistent by virtue of (\ref{eqn:coincidence}). 
The index $\nu_C(f)$ at a fixed curve $C \in X_1(f)$ is 
defined by 
\begin{equation} \label{eqn:nuC} 
\nu_C(f) := \nu_{\p}(f_x^*) 
\end{equation}
with a (any) point $x \in C \setminus I(f)$ and an (any) prime 
ideal $\p \in \A(C_x)$. 
\end{definition} 
\begin{remark} \label{rem:zero}
We have $\nu_x(f) > 0$ for at most finitely many points 
$x \in X_0(f)$. 
\end{remark}
\begin{definition} \label{def:type2} 
A fixed curve $C \in X_1(f)$ is said to be 
{\sl of type I} or {\sl of type I\!I} relative to 
$f : X \to X$ according as the prime ideal $\p \in \A(C_x)$ is 
of type I or of type I\!I relative to $f_x^* : A_x \to A_x$ 
in the sense of Definition \ref{def:type}. 
This definition does not depend on the choice of the point 
$x \in C \setminus I(f)$ and the ideal $\p \in \A(C_x)$ 
(see \cite[page 1016]{Saito}). 
Let $X_I(f)$ and $X_{I\!I}(f)$ denote the set of fixed curves 
of types I and the set of fixed curves of type I\!I 
respectively. 
Then there exists the direct sum decomposition as in 
(\ref{eqn:divide}). 
\end{definition} 
\par 
The preparation of all terminology and concepts needed to 
formulate Theorem \ref{thm:ssaito} is now complete. 
The separation condition (\ref{eqn:separate}) means that 
if $x \in X$ is an indeterminacy point of $f$ then $x$ is a 
holomorphic point of $f^{-1}$ and vice versa. 
Thus, under condition (\ref{eqn:separate}), all the 
fixed points needed to validate formula (\ref{eqn:formula}) 
are captured by the union $X_0(f) := X_0^{\circ}(f) \cup 
X_0^{\circ}(f^{-1})$. 
This is why we make the definition (\ref{eqn:fixed}). 
The consistency (\ref{eqn:coincidence}) with respect to 
$f^{\pm1}$ is well understood by the symmetry of the graphs 
$\varGamma_{f^{\pm1}}$ of $f^{\pm1}$, that is, by the fact that 
one graph is the reflection of the other in the diagonal 
$\varDelta$ of $X \times X$. 
This symmetry indicates the naturality of definition (\ref{eqn:nux}) 
since $\nu_{A_x}((f^{\pm1})_x^*)$ represent the 
degrees of intersection between $\varGamma_{f^{\pm1}}$ and 
$\varDelta$ at the point $(x,x)$. 
It is in these settings that Theorem \ref{thm:ssaito} is valid. 
A remark is in order at this stage. 
\begin{remark} \label{rem:separate} 
In dynamical situations the fixed point formula (\ref{eqn:formula}) 
is to be applied to the iterates $f^n$ of a map $f$, so that 
the separation condition (\ref{eqn:separate}) 
should be replaced by its iterated version: 
\[
I(f^n) \cap I(f^{-n}) = \emptyset \qquad \mbox{for every} \quad 
n \in \N. 
\]
It is easy to see that this condition follows from the AS 
condition (\ref{eqn:AS}). 
So the AS birational maps constitute a nice class of maps 
to which the fixed point formula (\ref{eqn:formula}) can 
be applied dynamically. 
\end{remark} 
\section{Stability of Indices} \label{sec:st}
The goal of this section is to establish Theorem \ref{thm:main}. 
This boils down to showing the following theorem in the 
abstract ring-theoretical setting as in the first part of 
Section \ref{sec:formula}. 
\begin{theorem} \label{thm:stability} 
Let $A := \C[\![z_1,z_2]\!]$ and $\si : A \to A$ a nontrivial 
continuous endomorphism in the $\m$-adic topology. 
If $\A_{I\!I}(\si)$ is nonempty, then for any $n \in \N$, 
\begin{itemize}
\item $\A(\si^n) = \A(\si)$, $\A_I(\si^n) = \A_I(\si)$ 
and $\A_{I\!I}(\si^n) = \A_{I\!I}(\si)$$;$ 
\item $\delta(\si^n) = \delta(\si)$ and 
$\nu_A(\si^n) = \nu_A(\si)$$;$
\item $\nu_{\p}(\si^n) = \nu_{\p}(\si)$ and 
$\mu_{\p}(\si^n) = \mu_{\p}(\si)$ for any $\p \in \A(\si)$. 
\end{itemize}
\end{theorem}
\par 
In order to prove this theorem we need some preliminaries. 
As in (\ref{eqn:si}), for each $n \in \N$ the endomorphism 
$\si^n : A \to A$ can be expressed as 
\begin{equation} \label{eqn:expand2}
\left\{
\begin{array}{rcl}
\si^n(z_1) &=& z_1 + g_n \cdot h_{n1}, \\[2mm]
\si^n(z_2) &=& z_2 + g_n \cdot h_{n2},
\end{array}
\right. 
\end{equation}
for some elements $g_n$, $h_{n1}$, $h_{n2} \in A$, where 
$g_n$ is nonzero and $h_{n1}$, $h_{n2}$ are relatively prime. 
By definition $\ag(\si^n) := (g_n)$ and $\bg(\si^n) := 
(h_{n1}, \, h_{n2})$ are the ideals generated by $g_n$ and by 
$h_{n1}$, $h_{n2}$, respectively. 
To simplify the notation we put $g := g_1$ and $h_i := h_{1i}$. 
\begin{lemma} \label{lem:dec} 
If $\A_{I\!I}(\si)$ is nonempty, then for any $n \in \N$ and 
$\p \in \A_{I\!I}(\si)$ we have: 
\[
\ag(\si^n) = \ag(\si), \qquad \bg(\si^n) = \bg(\si), \qquad 
h_{ni} \in n \cdot h_i + \p \cdot \bg(\si). 
\]
\end{lemma}
\proof 
We prove the lemma by induction on $n \in \N$. 
It is trivial for $n = 1$. 
Assume that the lemma holds for $n \in \N$. 
In what follows, to make the presentation simpler, we use the 
symbol $b_i$ to denote various elements of the ideal 
$\bg(\si) = (h_1, \, h_2)$. 
This abuse of notation causes no confusion when we are only 
interested in the argument modulo $\bg(\si)$. 
Let $g_{z_i}$ denote the formal partial derivative of $g$ with 
respect to $z_i$ $(i = 1,2)$. 
Considering the formal Taylor expansion of $\si^{n+1}(z_i) = 
\si(\si^n(z_i))$ with (\ref{eqn:si}) and (\ref{eqn:expand2}) 
taken into account, we have 
\begin{eqnarray}
\si^{n+1}(z_i) 
&=& \si(\si^n(z_i))  \nonumber \\[2mm]
&=& z_i + g \cdot h_i + 
g_n(z_1 + g \cdot h_1, \, z_2 + g \cdot h_2) \cdot 
h_{ni}(z_1 + g \cdot h_1, \, z_2 + g \cdot h_2) 
\nonumber \\[2mm]
&=& z_i + g \cdot h_i + 
g(z_1 + g \cdot h_1, \, z_2 + g \cdot h_2) \cdot 
h_{ni}(z_1 + g \cdot h_1, \, z_2 + g \cdot h_2) 
\nonumber \\[2mm]
&=& z_i + g \cdot h_i + ( g + g_{z_1} \cdot g \cdot h_1 + 
g_{z_2} \cdot g \cdot h_2 +g^2 \cdot b_i) \cdot 
(h_{ni} + g \cdot b_i) \nonumber \\[2mm]
&=& z_i + g \cdot \{ h_i + h_{ni} \cdot (1+ g_{z_1} \cdot h_1 
+ g_{z_2} \cdot h_2) + g \cdot b_i \} \nonumber \\[2mm] 
&=& z_i + g \cdot \wt{h}_{ni} \qquad (i = 1,2), 
\label{eqn:expand3} 
\end{eqnarray}
where we use the induction hypothesis $\ag(\si^n) = \ag(\si)$ 
in the third line and we put 
\begin{equation} \label{eqn:wtgni}
\wt{h}_{ni} := h_i + h_{ni} \cdot (1+ g_{z_1} \cdot h_1 
+ g_{z_2} \cdot h_2) + g \cdot b_i. 
\end{equation}
in the last line. 
We investigate the term $g_{z_1} \cdot h_1 + g_{z_2} \cdot h_2$. 
Take a prime element $p \in A$ such that $\p = (p)$. 
Since $\p = (p)$ divides $\ag(\si) = (g)$, the formal power 
series $g$ can be written 
\begin{equation} \label{eqn:GCD}
g = p^m \cdot \wt{g}
\end{equation}
for some $m \in \N$ and $\wt{g} \in A$. 
Since $\p \in \A_{I\!I}(\si)$, formula (\ref{eqn:a}) implies that 
the formal differential $\vp_{\si}$ in (\ref{eqn:omega}) is 
expressed as 
$\vp_{\si} := h_2 \cdot dz_1 - h_1 \cdot dz_2 = a \cdot d p + p 
\cdot (q_2 \cdot dz_1 - q_1 \cdot dz_2)$ with some elements 
$a$, $q_1$, $q_2 \in A$. 
Comparing the coefficients of $dz_1$ and $dz_2$ yields 
\begin{equation} \label{eqn:irr}
\left\{
\begin{array}{l}
h_1 = - a \cdot p_{z_2} + p \cdot q_1, \\[2mm]
h_2 = \phantom{-} a \cdot p_{z_1} + p \cdot q_2. 
\end{array}
\right.
\end{equation}
where $p_{z_i}$ is the formal partial derivative of $p$ with 
respect to $z_i$. 
By (\ref{eqn:GCD}) and (\ref{eqn:irr}), 
\[
\begin{array}{rcl}
g_{z_1} \cdot h_1+ g_{z_2} \cdot h_2 
&=& g_{z_1} \cdot ( - a \cdot p_{z_2} + p \cdot q_1 ) + 
g_{z_2} \cdot ( a \cdot p_{z_1} + p \cdot q_2 ) \\[2mm]
&=& ( m \cdot p^{m-1} \cdot p_{z_1} \cdot \wt{g} 
+ p^m \cdot \wt{g}_{z_1} ) 
\cdot (- a \cdot p_{z_2}) \\[2mm]
& & \phantom{ \{ } + ( m \cdot p^{m-1} \cdot p_{z_2} \cdot 
\wt{g} + p^m \cdot \wt{g}_{z_2} ) \cdot ( a \cdot p_{z_1} ) 
+ p \cdot (g_{z_1} \cdot q_1 + g_{z_2} \cdot q_2) \\[2mm] 
&=& a \cdot p^m \cdot (p_{z_1} \cdot \wt{g}_{z_2} - 
p_{z_2} \cdot \wt{g}_{z_1}) + 
p \cdot (g_{z_1} \cdot q_1 + g_{z_2} \cdot q_2) \\[2mm]
&=& p \cdot \{ a \cdot p^{m-1} \cdot (p_{z_1} \cdot \wt{g}_{z_2} - 
p_{z_2} \cdot \wt{g}_{z_1}) + 
(g_{z_1} \cdot q_1 + g_{z_2} \cdot q_2) \} \in \p. 
\end{array}
\]
Using this, $g \in \p$ and the induction hypothesis 
$h_{ni} \in n \cdot h_i + \p \cdot \bg(\si)$ in 
(\ref{eqn:wtgni}), we have
\begin{equation} \label{eqn:inc} 
\wt{h}_{ni}(z) \in (n+1) \cdot h_i(z) + \p \cdot \bg(\si) 
\qquad (i = 1, 2). 
\end{equation}
\par 
Now consider the ideal $\wt{\bg}(\si^n) := 
(\wt{h}_{n1}, \, \wt{h}_{n2})$ generated by $\wt{h}_{n1}$ and 
$\wt{h}_{n2}$. 
We show that 
\begin{equation} \label{eqn:ideal}
\wt{\bg}(\si^n) = \bg(\si).  
\end{equation}
The inclusion $\wt{\bg}(\si^n) \subset \bg(\si)$ is obvious, 
since $\wt{h}_{n1}$, $\wt{h}_{n2} \in \bg(\si)$ by (\ref{eqn:inc}). 
On the other hand, formula (\ref{eqn:inc}) also implies that 
there exist elements $r_{nij} \in \p$ $(i,j = 1, 2)$ such that 
\[
\left\{
\begin{array}{rclcl}
\wt{h}_{n1} 
&=& (n+1) \cdot h_1 + r_{n11} \cdot h_1 + r_{n12} \cdot h_2 
&=& (n + 1 + r_{n11}) \cdot h_1 + r_{n12} \cdot h_2, \\[2mm]
\wt{h}_{n2}
&=& (n+1) \cdot h_2 + r_{n21} \cdot h_1 + r_{n22} \cdot h_2 
&=& r_{n21} \cdot h_1 + (n + 1 + r_{n22}) \cdot h_2, 
\end{array}
\right. 
\]
If we put 
$r := (n+1+r_{n11}) \cdot (n+1+r_{n22})-r_{n12} \cdot r_{n21}$, 
then these equations yield 
\[
\left\{
\begin{array}{rcl}
r \cdot h_1 &=& 
(n+1+r_{n22}) \cdot \wt{h}_{n1} - r_{n12} \cdot \wt{h}_{n2}, 
\\[2mm]
r \cdot h_2 &=& 
r_{n21} \cdot \wt{h}_{n1} - (n+1+r_{n11}) \cdot \wt{h}_{n2}. 
\end{array}
\right.
\]
Since $r_{nij} \in \m$ $(i,j = 1, 2)$, the factor $r$ is an 
invertible element of $A$, so that one has 
$h_1$, $h_2 \in \wt{\bg}(\si^n)$. 
This yields the reverse inclusion $\bg(\si) \subset 
\wt{\bg}(\si^n)$ and the claim (\ref{eqn:ideal}) is proved. 
Since $h_1$ and $h_2$ are relatively prime, the equality 
(\ref{eqn:ideal}) implies that $\wt{h}_{n1}$ and $\wt{h}_{n2}$ 
are also relatively prime, so that from (\ref{eqn:expand3}) 
one can conclude that 
\begin{equation} \label{eqn:gn+1}
g_{n+1} = g, \qquad h_{n+1,i} = \wt{h}_{ni} \qquad (i =1, 2). 
\end{equation}
The first equality of (\ref{eqn:gn+1}) yields 
$\ag(\si^{n+1}) = \ag(\si)$. 
The second equality of (\ref{eqn:gn+1}) and (\ref{eqn:ideal}) 
lead to $\bg(\si^{n+1}) = \bg(\si)$. 
Finally the second equality of (\ref{eqn:gn+1}) and 
(\ref{eqn:inc}) give 
$h_{n+1,i} \in (n+1) \cdot h_i(z) + \p \cdot \bg(\si)$. 
Thus the induction is complete. \qed
\begin{lemma} \label{lem:typeI} 
If $\A_{I\!I}(\si)$ is nonempty then for any $\q \in \A_I(\si)$ 
there is $c_n \in \m$ such that 
\[
h_{ni} \in (n + c_n) \cdot h_i + \q. 
\]
\end{lemma}
\proof 
We prove the lemma by induction on $n \in \N$. 
Let $\p \in \A_{I\!I}(\si)$ be a prime ideal of type I\!I 
as in Lemma \ref{lem:dec} and its proof. 
It follows from (\ref{eqn:wtgni}) and (\ref{eqn:gn+1}) that 
\begin{equation} \label{eqn:rel}
h_{n+1,i} = h_i + h_{ni} \cdot (1+ g_{z_1} \cdot h_1 + g_{z_2} 
\cdot h_2) + g \cdot b_i \qquad (i = 1, 2), 
\end{equation}
where $b_i \in \bg(f^*)$. 
Let $q \in \m$ be a prime element such that $\q = (q)$. 
Since $\q$ is different from $\p$, the product $p \cdot q$ divides $g$. 
Moreover, since $p$, $q \in \m$, we have $g \in \q$ and $g \in \m^2$. 
If the assertion holds for $n$, then it follows from (\ref{eqn:rel}) 
that 
\[
\begin{array}{rcl}
h_{n+1,i} 
&=& h_i + h_{ni} \cdot (1+ g_{z_1} \cdot h_1 
+ g_{z_2} \cdot h_2) + g \cdot b_i \\[2mm]
&\equiv& h_i + h_{ni} \cdot (1+ g_{z_1} \cdot h_1 
+ g_{z_2} \cdot h_2) \qquad (\mathrm{mod} \,\, \q) \\[2mm]
&\equiv& h_i + (n + c_n) \cdot h_i \cdot (1+ g_{z_1} \cdot h_1 
+ g_{z_2} \cdot h_2) \qquad (\mathrm{mod} \,\, \q) \\[2mm]
&=& h_i \cdot \{ (n+1) + c_n + (n + c_n) \cdot (g_{z_1} \cdot h_1 
+ g_{z_2} \cdot h_2) \} \\[2mm]
&=& h_i \cdot \{ (n+1) + c_{n+1} \}, 
\end{array}
\]
where we use $g \in \q$ in the second line and the 
induction hypothesis in the third line, and we put 
$c_{n+1} := c_n + (n + c_n) \cdot 
(g_{z_1} \cdot h_1 + g_{z_2} \cdot h_2)$ in the last line. 
Since $g \in \m^2$, we have $g_{z_i} \in \m$ and hence 
$c_{n+1} \in \m$. 
Thus the assertion is true for $n+1$ and the 
induction is complete. \qed\par\medskip\noindent
{\it Proof of Theorem $\ref{thm:stability}$}. 
In view of (\ref{eqn:delta}), (\ref{eqn:nup}) and the definition 
of $\A(\si)$, the equalities $\ag(\si^n) = \ag(\si)$ and 
$\bg(\si^n) = \bg(\si)$ in Lemma \ref{lem:dec} imply that 
$\A(\si^n) = \A(\si)$, $\delta(\si^n) = \delta(\si)$ and 
$\nu_{\p}(\si^n) = \nu_{\p}(\si)$ for any $\p \in \A(\si)$. 
If $\p \in \A_{I\!I}(\si)$, then (\ref{eqn:omega}) and the last 
formula of Lemma \ref{lem:dec} imply that
\[
\vp_{\si^n} = h_{n2} \cdot dz_1 - h_{n1} \cdot dz_2 \in  
n \cdot (h_2 \cdot dz_1 - h_1 \cdot dz_2) + 
\p \cdot \hat{\Omega}_{A/\C}^1 
= n \cdot \vp_{\si} + \p \cdot \hat{\Omega}_{A/\C}^1, 
\]
and so $\tau_{\p}(\vp_{\si^n}) = n \cdot \tau_{\p}(\vp_{\si}) = 0$, 
which means that $\p \in \A_{I\!I}(\si^n)$. 
Moreover (\ref{eqn:a}) and (\ref{eqn:mup}) yield 
$\mu_{\p}(\si^n) = \mu_{\p}(\si)$. 
Next assume that $\p \in \A_I(\si)$ and rewrite $\p = \q$. 
By (\ref{eqn:omega}) and Lemma \ref{lem:typeI}, 
\[
\begin{array}{rcl}
\vp_{\si^n} = h_{n2} \cdot dz_1 - h_{n1} \cdot dz_2 
&\in& (n + c_n) \cdot (h_2 \cdot dz_1 - h_1 \cdot dz_2) 
+ \q \cdot \hat{\Omega}_{A/\C}^1 \\[2mm]
&=& (n + c_n) \cdot \vp_{\si} + \q \cdot \hat{\Omega}_{A/\C}^1, 
\end{array}
\]
where $n + c_n$ is an invertible element of $A$. 
Hence $\tau_{\q}(\vp_{\si^n}) = (n + c_n) \cdot 
\tau_{\q}(\vp_{\si}) \neq 0$, which means that 
$\q \in \A_I(\si^n)$. 
Moreover (\ref{eqn:a}) and (\ref{eqn:mup}) yield 
$\mu_{\q}(\si^n) = \mu_{\q}(\si)$. 
Therefore we have $\A_I(\si^n) = \A_I(\si)$, 
$\A_{I\!I}(\si^n) = \A_{I\!I}(\si)$ and 
$\mu_{\p}(\si^n) = \mu_{\p}(\si)$ for any $\p \in \A(\si)$. 
Finally the equality $\nu_A(\si^n) = \nu_A(\si)$ readily 
follows from (\ref{eqn:nuA}). 
The proof is complete. \qed \par\medskip 
We are now in a position to establish Theorem \ref{thm:main}. 
\par\medskip\noindent
{\it Proof of Theorem $\ref{thm:main}$}. 
Let $C \in X_{I\!I}(f)$ and take a point $x \in C \sm I(f)$. 
Then for each $n \in \N$ the map $f^n$ induces an endomorphism 
$(f^n)_x^* = (f_x^*)^n : A_x \to A_x$. 
Take any prime ideal $\p \in \A(C_x)$. 
Since $C \in X_{I\!I}(f)$, we have $\p \in \A_{I\!I}(f_x^*)$. 
Hence it follows from (\ref{eqn:nuC}) and 
Theorem \ref{thm:stability} that
\[
\nu_C(f^n) = \nu_{\p}((f_x^*)^n) = \nu_{\p}(f_x^*) = \nu_C(f), 
\]
which proves the first assertion of the theorem. 
Next, let $x \in X_0(f)$ be a fixed point of $f$ through which 
at least one fixed curve, say, $C \in X_{I\!I}(X)$ of type I\!I 
passes. 
In view of (\ref{eqn:fixed}) we may assume that 
$x \in X_0^{\circ}(f)$, namely, that $x \in C \setminus I(f)$; 
for, otherwise, we can replace $f$ by $f^{-1}$ and proceed 
in a similar manner. 
Now the endomorphisms $(f^n)_x^* = (f_x^*)^n : A_x \to A_x$ 
make sense and $\A_{I\!I}(f_x^*)$ is nonempty. 
Hence (\ref{eqn:nux}) and Theorem \ref{thm:stability} imply 
that 
\[
\nu_x(f^n) = \nu_{A_x}((f_x^*)^n) =\nu_{A_x}(f_x^*) = \nu_x(f), 
\]
which proves the second assertion of the theorem. 
Therefore Theorem \ref{thm:main} is established. 
\qed \par\medskip 
The following two remarks show that it is essential to assume that 
$x$ lies on a fixed curve of type I\!I in Theorem \ref{thm:main}. 
\begin{remark} \label{rem:ind1}
If $x \in X$ is an {\sl isolated} fixed point of all iterates 
$f^n$, then the indices $\nu_x(f^n)$ may depend on $n \in \N$. 
For example, consider a birational map $f : \P^2 \to \P^2$ 
expressed as 
\[
f(z_1, \, z_2) = (-2z_1-z_1^2-z_2, \, z_1), 
\]
in affine coordinates. 
Then the origin $(0,0)$ is an isolated fixed point of all 
iterates $f^n$. 
Indeed, assume the contrary that an iterate $f^n$ fixes some 
curve $C \subset \P^2$ passing through $(0,0)$. 
Let $L \subset \P^2$ be the line at infinity. 
We observe that $f$ has a superattracting fixed point 
$p^+ \in L$ and $f$ contracts $L$ into $p^+$. 
Similarly $f^{-1}$ has a superattracting fixed point $p^- \in L$ 
and $f^{-1}$ contracts $L$ into $p^-$, where $p^{\pm}$ are distinct. 
The curve $C$ intersects the line $L$ in a point, say, $q \in 
C \cap L$. 
If $q \neq p^-$ then $q = f^n(q) = p^+$, and if $q \neq p^+$ 
then $q = f^{-n}(q) = p^-$. 
But both equalities are impossible, because $p^+$ is an isolated 
fixed point of $f^n$ and $p^-$ is an isolated fixed 
point of $f^{-n}$. 
Thus $(0,0)$ is an isolated fixed point of all iterates $f^n$. 
A little calculation shows that $\nu_{(0,0)}(f) = 1$ and 
$\nu_{(0,0)}(f^2) = 3$ are distinct, though $\nu_{(0,0)}(f^n)$ 
are bounded by Theorem \ref{thm:SS}. 
Note that $f$ preserves the standard area form $dz_1 \wedge dz_2$. 
\end{remark}
\begin{remark} \label{rem:ind2}
If a point $x \in X$ is on a fixed curve $C$ of type I, 
then the indices $\nu_x(f^n)$ may depend on $n \in \N$. 
For example, consider a birational map $f : \P^2 \to \P^2$ 
expressed as 
\[
f(z_1, \, z_2) = (z_1+z_1(z_1^2+z_2), \, z_2+z_1^2). 
\]
in affine coordinates. 
Then $C := \{ z_1 = 0 \}$ is a fixed curve of type I. 
We can easily check that the index $\nu_{(0,z_2)} (f^n)$ at 
$(0,z_2) \in C$ is positive if and only if $z_2$ is a root of 
the equation $g_n(z_2) := (z_2+1)^n-1 = 0$. 
For example, $\nu_{(0,-2)}(f) = 0$ and $\nu_{(0,-2)}(f^2) \ge 1$ 
are distinct. 
Moreover, since the equation $g_n(z_2) = 0$ has $n$ distinct 
roots for each $n \in \N$, the number of points $x \in C$ 
such that $\nu_x(f) \ge 1$ grows linearly as $n$ tends to infinity. 
Thus, 
\[
\sum_{x \in C} \nu_{x}(f^n) \to + \infty \qquad (n \to +\infty). 
\]
On the other hand, if $C$ is a fixed curve of type I\!I, 
Theorem \ref{thm:main} and Remark \ref{rem:zero} imply that 
\[
\sum_{x \in C} \nu_{x}(f^n) = \sum_{x \in C} \nu_{x}(f) < \infty 
\qquad (n \in \N). 
\]
In order to get the last equality we need the invariance of 
the indices $\nu_x(f^n)$ as a function of $n \in \N$, whereas 
the boundedness of $\nu_x(f^n)$ as in Theorem \ref{thm:SS} is 
not enough for this aim. 
\end{remark} 
\section{Finiteness of Periodic Curves} \label{sec:Fib}
The aim of this section is to discuss the finiteness of the number 
of periodic curves of type I\!I and especially to prove 
Theorem \ref{thm:main2}. 
Actually we establish more refined results 
(Theorems \ref{thm:pos-kodaira}, \ref{thm:rational} and 
\ref{thm:zero-deg}), dividing our discussion 
into three cases according to the values of the first dynamical 
degree $\lambda(f)$ of $f$ and also to the values of the 
Kodaira dimension $\mathrm{kod}(X)$ of $X$. 
Theorem \ref{thm:main2} is then deduced as a corollary of these 
results. 
We begin this section with a result on the constraints for the prime 
periods of two intersecting periodic curves. 
It plays an important role in the main discussion of this section, 
while it is also of intrinsic interest in its own light. 
\begin{theorem} \label{thm:prime} 
Let $f : X \to X$ be a nontrivial AS birational map and let 
$C$ be a periodic curve of type I\!I with prime period $n$. 
If $C'$ is a periodic curve of prime period $m$ that intersects 
$C$, then $m$ is a divisor of $n$. 
If moreover $C'$ is of type I\!I, then $m = n$. 
\end{theorem}
\proof
Assume that $C$ and $C'$ intersect in a point $x \in C \cap C'$. 
In view of (\ref{eqn:fixed}), since the map $f$ is assumed to be 
AS, some choice of double signs $(\ve, \delta) \in \{\pm1\}^2$ 
makes $x \in X_0^{\circ}(f^{\ve n}) \cap X_0^{\circ}(f^{\delta \m})$. 
Among the four cases we only discuss the two cases $(\ve, \delta) = 
(+,+)$ and $(\ve, \delta) = (+,-)$, as the remaining cases 
$(\ve, \delta) = (-,+)$ and $(\ve, \delta) = (-,-)$ can be 
treated in similar manners. \par 
First, assume that $x \in X_0^{\circ}(f^n) \cap 
X_0^{\circ}(f^m)$, namely, that $x \in C \sm I(f^n)$ and 
$x \in C' \sm I(f^m)$. 
Then one can think of three endomorphisms: 
\[
(f^n)_x^*, \quad (f^m)_x^*, \quad  
((f^n)_x^*)^m = (f^{nm})_x^* = ((f^m)_x^*)^n \, :\, A_x \to A_x. 
\]
Since $C \in X_{I\!I}(f^n)$ passes through $x$, any irreducible 
component of the germ $C_x$ at $x$ defines an element of 
$\A_{I\!I}((f^n)_x^*)$, which is therefore nonempty. 
So Theorem \ref{thm:stability} with $\si = (f^n)_x^*$ yields 
\[
\A((f^n)_x^*) = \A(((f^n)_x^*)^m) = \A(((f^m)_x^*)^n) 
\supset \A((f^m)_x^*). 
\]
The prime ideal $\p$ corresponding to any irreducible component 
of the germ $C_x'$ is an element of $\A((f^m)_x^*)$. 
Thus the inclusion relation above yields $\p \in \A((f^n)_x^*)$, 
which means that $C' \sm I(f^n)$ is fixed pointwise by $f^n$. 
Now recall that $C' \sm I(f^m)$ is fixed pointwise by $f^m$. 
Write $n = k m + r$ with $k \in \Z_{\ge 0}$ and 
$r \in \{0,1,\dots,m-1\}$. 
Then $C' \sm (I(f^n) \cup I(f^m) \cup I(f^r))$ and hence 
$C' \sm I(f^r)$ are fixed pointwise by $f^r$. 
Since $C'$ is a periodic curve of prime period $m > r$, 
we must have $r = 0$ and $n = km$. 
Hence $m$ is a divisor of $n$. \par 
Secondly, assume that $x \in X_0^{\circ}(f^n) \cap 
X_0^{\circ}(f^{-m})$, namely, that $x \in C \sm I(f^n)$ and 
$x \in C' \sm I(f^{-m})$. 
Then we have $x \in X_0^{\circ}(f^{nm}) \cap 
X_0^{\circ}(f^{-mn})$ and hence $f^{nm}$ defines a local 
biholomorphism around $x$, which induces a ring 
automorphism $((f^n)_x^*)^m = (f^{nm})_x^* : A_x \to A_x$ 
together with its inverse $((f^{-m})_x^*)^n = 
(f^{-nm})_x^* = ((f^{nm})_x^*)^{-1} : A_x \to A_x$. 
By Lemma \ref{lem:inverse} we have 
$\A((f^{nm})_x^*) = \A(((f^{nm})_x^*)^{-1})$ and hence 
$\A(((f^n)_x^*)^m) = \A(((f^{-m})_x^*)^n) \supset 
\A((f^{-m})_x^*)$. 
On the other hand, since $C \in X_{I\!I}(f^n)$ passes through 
$x$, the set $\A_{I\!I}((f^n)_x^*)$ is nonempty. 
So Theorem \ref{thm:stability} with $\si = (f^n)_x^*$ implies 
$\A((f^n)_x^*) = \A(((f^n)_x^*)^m) \supset \A((f^{-m})_x^*)$. 
Now note that $X_1(f^m) = X_1(f^{-m})$ by (\ref{eqn:symm}). 
Since $C' \in X_1(f^m) = X_1(f^{-m})$ passes through $x$, 
any irreducible component of the germ $C'_x$ defines a prime 
element $\p \in \A((f^{-m})_x^*)$. 
By the inclusion above we have $\p \in \A((f^n)_x^*)$, which 
means that $C' \sm I(f^n)$ is fixed pointwise by $f^n$, while 
$C' \sm I(f^m)$ is fixed pointwise by $f^m$. 
The remaining argument is the same as in the last paragraph. 
We have $n = km$ for some $k \in \Z_{\ge 0}$. \par 
In any case it is shown that $m$ is a divisor of $n$. 
If moreover $C'$ is of type I\!I, then the same reasoning as 
above with $C$ replaced by $C'$ implies that $n$ is a divisor 
of $m$ and hence $m = n$. \qed\par\medskip
Diller and Favre \cite{DF} give a classification of bimeromorphic 
maps on a compact K\"ahler surface in terms of their first dynamical 
degrees (\ref{eqn:FDD}). 
We make use of this classification in our discussion. 
\begin{theorem}[\cite{DF}] \label{thm:DF} 
Let $f : X \to X$ be a bimeromorphic map on a compact 
K\"{a}hler surface $X$. 
Then $f$ is classified in the following manner up to 
bimeromorphic conjugacy. 
\par\smallskip 
$\bullet$ If $\lambda(f) = 1$, then exactly one of the 
following is true: 
\begin{enumerate}
\item[$(0)$] $|\!|(f^n)^*|\!|$ are bounded and $f^m$ is an 
automorphism isotopic to the identity for some $m \in \N$$;$  
\item[$(1)$] $|\!|(f^n)^*|\!|$ grow linearly and $f$ preserves a 
unique rational fibration $\pi : X \to S$$;$   
\item[$(2)$] $|\!|(f^n)^*|\!|$ grow quadratically and $f$ 
preserves a unique elliptic fibration $\pi : X \to S$. 
\end{enumerate}
\par\smallskip 
$\bullet$ If $\lambda(f) > 1$, then either 
\begin{enumerate}
\item[$(3)$] $X$ is a rational surface with $f$ an automorphism 
or merely a bimeromorphic map$;$ or 
\item[$(4)$] $f$ is an automorphism of a $\mathrm{K3}$ surface, 
an Enriques surface or a complex $2$-torus. 
\end{enumerate}
\end{theorem} 
\par 
The case $\lambda(f) = 1$ of low dynamical complexity is more 
or less easy to handle. 
So we are mostly concerned with the case $\lambda(f) > 1$, 
which is divided into subcases (3) and (4) according to 
$\mathrm{kod}(X) = - \infty$ and $\mathrm{kod}(X) \ge 0$ 
respectively. 
We begin with the last case (4). 
\begin{theorem} \label{thm:pos-kodaira} 
If $X$ is a compact K\"{a}hler surface of Kodaira dimension 
$\mathrm{kod}(X) \ge 0$ and $f : X \to X$ is a bimeromorphic 
map of first dynamical degree $\lambda(f) > 1$, then $f$ has no 
irreducible periodic curves of nonnegative self-intersection. 
\end{theorem} 
\proof
It suffices to show that $f$ has no irreducible fixed curves of 
nonnegative self-intersection, because if one wants to consider 
periodic curves of period $n$, then one may replace $f$ with $f^n$ 
upon noting that $\lambda(f^n) = \lambda(f)^n > 1$. 
Now assume the contrary that $f$ admits an irreducible fixed 
curve $C$ of nonnegative self-intersection $C^2 \ge 0$. 
Since $\mathrm{kod}(X) \ge 0$, it follows from 
\cite[Chap. VI, (1.1) Theorem]{BHPV} and 
\cite[Proposition 7.5]{DF} that there exists a commutative 
diagram 
\[
\begin{CD}
X @> f >> X \\
@V \varphi VV  @VV \varphi V \\
Y @>> g > Y 
\end{CD}
\]
such that $\varphi$ is a proper modification, $Y$ is the unique 
minimal model of $X$ and $g$ is an automorphism such that 
$\lambda := \lambda(g) = \lambda(f) > 1$. 
Let $E(\varphi)$ be the exceptional set of $\varphi$. 
Since any irreducible component of $E(\varphi)$ 
has a negative self-intersection, the curve $C$ is not 
contained in $E(\varphi)$. 
So $C' := \varphi(C)$ is an irreducible fixed curve of $g$. 
Since blowing down a curve does not decrease its 
self-intersection number, 
$C'$ has also a nonnegative self-intersection. 
By \cite[Theorem 0.3]{DF} there exists a nef class 
$\theta \in H^{1,1}(Y)$ such that $g^* \theta = \lambda \theta$. 
Then we have $g^* C' = C'$ and hence 
$\lambda (C', \theta) = (C', \lambda \theta) = 
(g^*C', g^* \theta) = (C', \theta)$, where in the last equality 
we use the fact that an automorphism preserves the 
intersection form. 
But, since $\lambda > 1$, we have $(C', \theta) = 0$ together 
with $(C')^2 \ge 0$ and $\theta^2 \ge 0$. 
Then the Hodge index theorem and 
\cite[Chap. IV, (7.2) Corollary]{BHPV} imply that 
$(C')^2 = \theta^2 = (C', \theta) = 0$ and there exists a 
positive constant $a > 0$ such that 
$C' = a \, \theta$ in $H^{1,1}(Y)$. 
Applying $g^*$ to this equality, we have 
$C' = a \lambda \theta$. 
But, since $\lambda > 1$, we have $C' = \theta = 0$ in 
$H^{1,1}(Y)$. 
This is a contradiction. \qed 
\begin{remark} \label{rem:pos-kodaira} 
If $\mathrm{kod}(X) = - \infty$, Theorem \ref{thm:pos-kodaira} 
is not true in general. 
For example, consider the birational map $f : \P^2 \to \P^2$ 
discussed in Remark \ref{rem:ind2}. 
We can see that $f^{-m}I(f) = \{ [0:0:1] \}$ and 
$f^{n}I(f^{-1}) = \{ [0:1:0], \, [1:0:-1] \}$ for every 
$m$, $n \ge 0$ and so $f$ is AS. 
Since $\deg f = 3$, we have $\lambda(f) = 3 > 1$. 
The map $f$ has the line $C := \{z_1 = 0 \}$ as a fixed curve 
of self-intersection $C^2 = 1$. 
This gives a counterexample to Theorem \ref{thm:pos-kodaira} 
when $X$ is a rational surface, or more specifically when 
$X = \P^2$. 
We can also construct an AS birational map of a rational surface 
having a fixed curve of zero self-intersection by blowing up a 
suitable point of $C$, say, $[1:0:0] \in C$, and lifting the map 
$f$ to the surface upstairs. 
\end{remark}
\par 
We proceed to the case where $\lambda(f) > 1$ and 
$\mathrm{kod}(X) = -\infty$, namely, the case (3). 
\begin{theorem} \label{thm:rational}
If $X$ is a smooth rational surface and $f : X \to X$ is an 
AS birational map of first dynamical degree $\lambda(f) > 1$, 
then all the irreducible periodic curves of type I\!I of $f$ 
with zero self-intersection have one and the same prime period. 
\end{theorem}
\par 
The proof is divided into several steps and begins with some 
generality on fibrations. 
\begin{lemma} \label{lem:zari}
Let $\pi : X \to S$ be a fibration with connected fibers of a 
smooth surface $X$ and let $C \subset X$ be a curve with zero 
self-intersection such that $\pi(C)=\{t\}$. 
Then the following hold: 
\begin{enumerate}
\item $X_t := \pi^{-1}(t) = s \, C$ for some $s \in \N$.
\item If a connected curve $C'$ is disjoint from $C$,  
then $C'$ is contained in a fiber of $\pi$. 
\end{enumerate}
\end{lemma}
\proof
Assertion (1): Write $X_t = s \, C + D$, where $D$ is an effective 
divisor not containing $C$. 
Assume that $D$ is nonempty. 
Since each fiber of $\pi$ is connected, $C$ and $D$ intersect 
so that $C \cdot D > 0$. 
Then $(m C + D)^2 = 2m C \cdot D + D^2 > 0$ for a sufficiently 
large integer $m \in \N$. 
This contradicts the fact that any divisor supported on a fiber 
has a nonpositive self-intersection 
(see \cite[Chap. III, (8.2) Lemma]{BHPV}). 
Hence $D$ must be empty and assertion (1) is proved. \par 
Assertion (2): 
Since $C'$ is connected, its image $\pi(C')$ is a connected 
algebraic subset of $S$, which must be a single point of $S$ or 
the entire curve $S$. 
On the other hand, since $C'$ is disjoint from $C$, 
assertion (1) implies that $C'$ is also disjoint from $X_t$ 
and hence $\pi(C') \subset S \sm \{t\}$. 
Therefore $\pi(C')$ must be a single point, so that 
$C'$ is contained in a fiber of $\pi$. \qed\par\medskip  
To prove Theorem \ref{thm:rational} by contradiction, 
we assume the contrary and proceed as follows. 
\begin{lemma} \label{lem:fib1} 
Let $X$ be a smooth rational surface and $f : X \to X$ an AS 
birational map with $\lambda(f) > 1$. 
Assume that $f$ admits two irreducible periodic curves $C_1$ and 
$C_2$ of type I\!I with zero self-intersection and with distinct 
prime periods $n_1$ and $n_2$ respectively. 
Then there exists a fibration $\pi : X \to S$ with connected fibers 
such that $C_1$ and $C_2$ are fibers of $\pi$. 
Moreover, if a connected curve $C \subset X$ is disjoint from 
$C_1$ or $C_2$, then $C$ is contained in a fiber of $\pi$. 
\end{lemma}
\proof 
Since $C_1$ and $C_2$ are periodic curves of type I\!I with 
distinct prime periods, Theorem \ref{thm:prime} implies that 
$C_1 \cap C_2 = \emptyset$. 
This together with the assumption of zero self-intersection 
yields $C_1^2 = C_2^2 = (C_1, C_2) = 0$. 
Then Hodge index theorem tells us that $C_1$ and $C_2$ 
are linearly dependent in $\mathrm{NS}(X) \otimes \R$, where 
$\mathrm{NS}(X)$ is the N\'{e}ron-Severi group of $X$. 
Moreover, since $C_1$ and $C_2$ are positive and 
$\mathrm{NS}(X)$ is defined over $\Z$, there exist $a_1$, 
$a_2 \in \N$ such that $a_1 C_1 = a_2 C_2$ in $\mathrm{NS}(X)$. 
Since $X$ is a rational surface, we have 
$H^1(X,\mathcal{O}_X) = 0$, so that the first Chern class 
map $c_1 : \mathrm{Pic}(X) \cong H^1(X,\mathcal{O}_X^*) \to 
\mathrm{NS}(X)$ is injective. 
Thus $a_1 C_1 = a_2 C_2$ in $\mathrm{Pic}(X)$, namely, the 
divisors $a_1 C_1$ and $a_2 C_2$ are linearly equivalent. 
Then there exists a surjective holomorphic map 
$\wt{\pi} : X \to \P^1$ such that $\wt{\pi}^{-1}(0) = 
a_1C_1$ and $\wt{\pi}^{-1}(\infty) = a_2 C_2$. 
Using the Stein factorization 
(see \cite[Chap. I, (8.1) Theorem]{BHPV}), 
we obtain a fibration $\pi : X \to S$ over some curve $S$ with 
connected fibers and a finite morphism $\phi : S \to \P^1$ such 
that $\wt{\pi} = \phi \circ \pi : X \to \P^1$. 
Since $C_i$ is a connected fiber of $\wt{\pi}$, 
it is also a fiber of the fibration $\pi$. 
Finally, if a connected curve $C$ is disjoint from $C_i$ 
for at least one $i \in \{1, 2\}$, then assertion (2) of 
Lemma \ref{lem:zari} implies that $C$ is contained in a fiber 
of the fibration $\pi$. 
\qed \par\medskip
For any curve $C$ on a surface $X$, one can define the 
pushforward $f_*C$ of $C$ by pulling back the local defining 
function of $C$ by $f^{-1}$. 
In general, $f_*C - fC$ is a nonnegative linear combination of 
irreducible components of $E(f^{-1})$. 
In this situation, Diller and Favre \cite[Corollary 3.4]{DF} 
obtain a useful formula for the intersection number of 
the pushforwards of two curves: 
\begin{equation} \label{eqn:PPF} 
(f_*C, f_*C')=(C, C') + Q(C, C'), 
\end{equation}
where $Q(C, C')$ is a nonnegative Hermitian form expressed as 
\[
Q(C, C') = \sum_{\underset{\text{irreducible}}{V \subset E(f):}} 
k(V) \cdot (C, V) \cdot (C', V), 
\]
with some positive integer $k(V) \in \N$ for each irreducible 
component $V$ of $E(f)$. 
Note that $Q(C,C) = 0$ if and only if $(C, V) = 0$ for any 
irreducible component $V$ of $E(f)$. 
\begin{lemma} \label{lem:fib2}
Under the assumptions of Lemma $\ref{lem:fib1}$, for any 
$i \ge 0$, the $i$-th pushforward $f_*^iC_1$ has zero 
self-intersection$;$ $(f_*^iC_1, V) = 0$ for any irreducible 
component $V$ of $E(f)$$;$ and $f_*^iC_1$ is equal to a fiber 
of the fibration $\pi : X \to S$ constructed in 
Lemma $\ref{lem:fib1}$. 
\end{lemma}
\proof
For each $i \ge 0$ we put $f_*^iC_1 = f^iC_1 + E_i$ with some 
effective divisor $E_i$. 
Since $f^iC_1$ is a periodic curve of prime period $n_1$, 
Theorem \ref{thm:prime} implies that $f^iC_1$ is disjoint from 
$C_2$. 
We see that $E_i$ is also disjoint from $C_2$. 
Indeed, if $E_i$ intersects $C_2$, then $f^{-i}C_2$ must 
intersect $C_1$, but this contradicts the fact that $f^{-i}C_2$ 
is a periodic curve of prime period $n_2$. 
Thus $f_*^iC_1 = f^iC_1 + E_i$ is disjoint from $C_2$. 
Since $f_*^iC_i$ is connected, Lemma \ref{lem:fib1} implies 
that $f_*^iC_1$ is contained in a fiber of the fibration $\pi$. 
In particular, $f_*^{n_1}C_1 = f^{n_1}C_1 + E_{n_1} = C_1 + E_{n_1}$ 
is contained in a fiber of $\pi$. 
Because the fiber containing $C_1$ is $C_1$ itself by 
Lemma \ref{lem:fib1}, we have $E_{n_1} = \emptyset$ and thus 
$f_*^{n_1}(C_1) = C_1$. 
A repeated use of formula (\ref{eqn:PPF}) yields 
\[ 
\begin{array}{rcl}
C_1^2 &=& (f_*^{n_1}C_1)^2 = (f_*^{n_1-1}C_1)^2 + 
Q(f_*^{n_1-1}C_1,f_*^{n_1-1}C_1) = \cdots \\[2mm]
&=& (f_*^jC_1)^2 + 
\displaystyle \sum_{i=j}^{n_1-1} Q(f_*^{i}C_1,f_*^{i}C_1) = \cdots 
\\[3mm]
&=& C_1^2 + \displaystyle \sum_{i=0}^{n_1-1} Q(f_*^{i}C_1,f_*^{i}C_1).  
\end{array}
\]
Since $Q$ is a nonnegative Hermitian form, this formula shows that 
$Q(f_*^{i}C_1,\, f_*^{i}C_1) = 0$ for any $0 \le i \le n_1-1$ and 
thus $(f_*^iC_1)^2 = 0$ and $(f_*^iC_1, V) = 0$ for any 
$0 \le i \le n_1-1$ and any irreducible component $V$ of $E(f)$. 
These are true for any $i \ge 0$ since $f_*^{n_1}C_1 = C_1$. 
By assertion (1) of Lemma \ref{lem:zari} we can conclude that 
$f_*^iC_1$ is a fiber of $\pi$. 
The proof is complete. \qed
\begin{lemma} \label{lem:fib3}
Under the assumptions of Lemma $\ref{lem:fib1}$, the map $f : X \to X$ 
preserves the fibration $\pi : X \to S$ constructed in 
Lemma $\ref{lem:fib1}$. 
\end{lemma}
\proof
Fix a point $t \in S$ and consider the fiber $X_t:=\pi^{-1}(t)$ over $t$. 
We show that $f_* X_t$ is contained in a fiber of the fibration $\pi$. 
This is true for $X_t = C_1$ by Lemma \ref{lem:fib1}, so that 
we may assume that $X_t \cap C_1 = \emptyset$. 
By formula (\ref{eqn:PPF}) and Lemma \ref{lem:fib2}, we have 
\[
(f_*X_t, f_*C_1) = (X_t, C_1) + 
\sum_{\underset{\text{irreducible}}{V \subset E(f):}} 
k(V) \cdot (X_t, V) \cdot (C_1,V) = (X_t,C_1) = 0. 
\]
This means that $f_*X_t$ is disjoint from $f_*C_1$. 
Since $f_*C_1$ is equal to a fiber of $\pi$ by Lemma \ref{lem:fib2}, 
$f_*X_t$ is contained in a fiber of $\pi$ by Lemma \ref{lem:zari}. 
\qed \par\medskip\noindent
{\it Proof of Theorem $\ref{thm:rational}$}. 
Assume the contrary that the theorem does not hold. 
Then the assumptions of Lemma \ref{lem:fib1} are satisfied 
and hence $f$ preserves a fibration $\pi : X \to S$ by Lemma 
\ref{lem:fib3}. 
But this is absurd, because a bimeromorphic map preserving 
a fibration has first dynamical degree $\lambda(f) = 1$ 
(see \cite[Corollary 1.3]{CF}). 
Thus the theorem is established. \qed
\begin{theorem} \label{thm:zero-deg}
Let $X$ be a compact K\"{a}hler surface and $f : X \to X$ a 
bimeromorphic map with $\lambda(f) = 1$ such that $f^n$ is 
not isotopic to the identity for any $n \in \N$. 
Let $\pi : X \to S$ be the unique rational or elliptic 
fibration preserved by $f$ in Theorem $\ref{thm:DF}$. 
If $f$ admits two irreducible periodic curves $C_1$ and $C_2$ of 
type I\!I with zero self-intersection and with distinct prime 
periods $n_1$ and $n_2$ respectively, then any irreducible 
periodic curve of type I\!I with an arbitrary prime period and 
any irreducible periodic curve of type I whose prime period is 
not a common divisor of $n_1$ and $n_2$ are contained in fibers 
of the fibration $\pi$. 
\end{theorem}
\proof 
First notice that $C_1$ and $C_2$ are disjoint by Theorem 
\ref{thm:prime}, since $n_1$ and $n_2$ are distinct. 
Let $C$ be a periodic curve of type I\!I with an arbitrary prime 
period $n$ or a periodic curve of type I whose prime period $n$ 
is not a common divisor of $n_1$ and $n_2$. 
We claim that $C$ is disjoint form either $C_1$ or $C_2$. 
Indeed, if $C$ is of type I\!I, then its prime period $n$ 
is different from $n_i$ for at least one $i \in \{1,2\}$, 
and hence Theorem \ref{thm:prime} shows that $C$ is 
disjoint from $C_i$. 
On the other hand, if $C$ is of type I and meets both $C_1$ and 
$C_2$, then $n$ must divide both $n_1$ and $n_2$, which 
contradicts the assumption by Theorem \ref{thm:prime}. 
Hence the claim is verified. \par 
We establish the theorem by a case-by-case check. 
First, when $\pi : X \to S$ is an elliptic fibration, it 
follows from \cite[Theorem 3.4]{DJS} that $C$ is contained 
in a fiber of the fibration $\pi$. 
Next we consider the case where $\pi : X \to S$ is a rational 
fibration. 
Assume that $\pi(C_i) = S $ and $C_j$ is contained in a 
fiber of $\pi$ for some $\{i,j\}=\{1, 2\}$, then $C_j$ is equal 
to a fiber of $\pi$ by assertion (1) of Lemma \ref{lem:zari} and 
thus $C_i$ must intersect $C_j$, but this contradicts the fact 
that $C_1 \cap C_2 = \emptyset$. 
Now assume that $\pi(C_1) = \pi(C_2) = S$. 
Then any fiber $X_t := \pi^{-1}(t)$ of $\pi$ meets both 
$C_1$ and $C_2$. 
If we take $t \in S$ to be sufficiently generic, then 
$X_t \cong \P^1$ and an intersection point 
$p_i \in X_t \cap C_i$ becomes a periodic point of $f$ with 
prime period $n_i$ for each $i \in \{1, 2\}$. 
Since $f$ preserves the fibration $\pi$ and 
$f^{n_i}(p_i) = p_i$, we have two automorphisms 
$f^{n_i}|_{X_t} : X_t \to X_t$ $(i = 1, 2)$. 
If $d$ denotes the greatest common divisor of $n_1$ and $n_2$, 
then $f^{d}|_{X_t}$ becomes an automorphism of $X_t$ having $p_1$ 
and $p_2$ as periodic points of prime periods $n_1/d$ and $n_2/d$ 
respectively. 
Hence $f^{d}|_{X_t}$ is a linear fractional transformation with 
two periodic points of distinct prime periods, 
but this is impossible. 
Thus for each $i = 1$, $2$, the curve $C_i$ must be contained in 
a fiber of $\pi$ and in fact equal to that fiber by assertion (1) 
of Lemma \ref{lem:zari}. 
Since $C$ is disjoint from either $C_1$ or $C_2$, assertion (2) 
of Lemma \ref{lem:zari} implies that $C$ is contained in 
a fiber of $\pi$.  \qed \par\medskip
Finally, in order to prove Theorem \ref{thm:main2}, we need the 
following lemma. 
\begin{lemma} \label{lem:picard}
Let $k \in \N$. 
Given $\rho(X) + k$ irreducible periodic curves of type I\!I 
with mutually distinct prime periods, then at least $k$ of 
them have zero self-intersection. 
\end{lemma}
\proof 
Let $C_1, \dots, C_{\rho(X)+k}$ be the periodic curves 
of type I\!I with mutually distinct prime periods. 
They are mutually disjoint by Theorem \ref{thm:prime}. 
If the contrary to the lemma holds, then we may assume that 
$C_i^2$ is nonzero for every $i \in \{1, 2, \cdots, \rho(X)+1\}$ 
after rearranging the suffixes if necessary. 
Since $\rho(X) = \mathrm{dim}_{\R} \, \mathrm{NS}(X) \otimes \R$, 
there is a nontrivial linear relation 
$r_1 C_1+ \cdots + r_{\rho(X)+1} C_{\rho(X)+1} = 0$ in 
$\mathrm{NS}(X) \otimes \R$. 
Since $C_i$ and $C_j$ are disjoint for every distinct 
$i$ and $j$, the linear relation yields 
\[
0 = \sum_{j = 1}^{\rho(X)+1} r_j (C_i, C_j) = r_i C_i^2 
\qquad (i = 1, \dots, \rho(X)+1), 
\]
which means that $r_1 = \cdots = r_{\rho(X)+1} = 0$. 
This contradicts the fact that the linear relation is nontrivial. 
Thus at least $k$ of $C_1, \dots, C_{\rho(X)+k}$ have zero 
self-intersection. \qed\par\medskip\noindent
{\it Proof of Theorem $\ref{thm:main2}$}. 
Assertion (1). 
First we consider the case $\mathrm{kod}(X) \ge 0$. 
Assume that $f$ admits $\rho(X) + 1$ irreducible periodic curves of 
type I\!I with mutually distinct prime periods. 
Then Lemma \ref{lem:picard} with $k = 1$ implies that at least one 
of them has zero self-intersection. 
But this is impossible by Theorem \ref{thm:pos-kodaira}. 
Thus there are at most $\rho(X)$ irreducible periodic curves of 
type I\!I with mutually distinct prime periods. 
This proves the item (2) of Remark \ref{rem:main2}. 
Next we consider the case $\mathrm{kod}(X) = - \infty$, namely, the 
case where $X$ is rational. 
Assume that $f$ admits $\rho(X) + 2$ irreducible periodic curves of 
type I\!I with mutually distinct prime periods. 
Then Lemma \ref{lem:picard} with $k = 2$ implies that at least two 
of them have zero self-intersection. 
But this is impossible by Theorem \ref{thm:rational}. 
Thus there are at most $\rho(X)+1$ irreducible periodic curves of 
type I\!I with mutually distinct prime periods. 
Therefore assertion (1) of the theorem is proved. 
\par 
Assertion (2). 
Assume that $f$ has more than $\rho(X) + 1$ irreducible periodic 
curves of type I\!I with mutually distinct prime periods. 
Again by Lemma \ref{lem:picard} with $k = 2$, at least two of 
them have zero self-intersection. 
Then Theorem \ref{thm:zero-deg} implies that any irreducible 
periodic curve of type I\!I is contained in a fiber of the 
fibration $\pi$. \qed 
\section{Area-Preserving Maps} \label{sec:APD}
The aim of this section is to discuss the absence of periodic 
curves of type I for an area-preserving map and to 
prove Theorems \ref{thm:main2}. 
Given a fixed curve $C \in X_1(f)$, we take a smooth 
point $x$ of $C$ and identify $A_x$ with $\C[\![z_1,z_2]\!]$ in 
such a manner that $C$ has the local defining equation $z_1 = 0$ 
near $x$. 
Then the induced endomorphism $f_x^* : A_x \to A_x$ can be 
expressed as 
\begin{equation} \label{eqn:expand}
\left\{
\begin{array}{rcl}
f_x^*(z_1) &=& z_1 + z_1^k \cdot f_1, \\[2mm]
f_x^*(z_2) &=& z_2 + z_1^l \cdot f_2,
\end{array}
\right. 
\end{equation}
for some $k$, $l \in \N \cup \{\infty\}$ and some $f_i \in A_x$ 
such that $f_i(0,z_2)$ is a nonzero element of $\C[\![z_2]\!]$. 
Here we put $z_1^{\infty} := 0$ by convention and we remark 
that at least one of $k$ and $l$ is finite. 
\begin{lemma} \label{lem:typeII}
For $C \in X_1(f)$, we have $\nu_{C}(f) = \min \{ k,l \}$ and 
$C \in X_{I\!I}(f)$ if and only if $k > l$. 
\end{lemma}
\proof
In (\ref{eqn:expand}) we put $f_i = u \cdot v_i$, where $v_1$ and 
$v_2$ are relatively prime. 
If $k > l$, then comparing (\ref{eqn:si}) with (\ref{eqn:expand}) 
we have $g = z_1^l \cdot u$, $h_1 = z_1^{k-l} \cdot v_1$ and 
$h_2 = v_2$. 
Hence (\ref{eqn:nup}) and (\ref{eqn:omega}) yield 
\begin{eqnarray*}
\nu_{C}(f) 
&=& \max \{\,m \in \N \, | \, (g) \subset (z_1)^m \,\} = l, 
\\[2mm] 
\vp_{f_x^*} 
&=& v_2 \cdot dz_1 - z_1^{k-l} \cdot v_1 \cdot dz_2 \in 
\hat{\Omega}_{A_x/\C}^1, 
\end{eqnarray*}
which shows that $\tau_{(z_1)}(\vp_{f_x^*}) = 0$ and hence $C$ is 
of type I\!I relative to $f$. 
On the other hand, if $k \leq l$, then we have $g = z_1^k \cdot u$, 
$h_1 = v_1$ and $h_2 = z_1^{l-k} \cdot v_2$. 
Hence (\ref{eqn:nup}) and (\ref{eqn:omega}) yield 
\begin{eqnarray*}
\nu_{C}(f) 
&=& \max \{\, m \in \N \, | \, (g) \subset (z_1)^m \,\} = k, 
\\[2mm] 
\vp_{f_x^*} 
&=& z_1^{l-k} \cdot v_2 \cdot dz_1 - v_1 \cdot dz_2 
\in \hat{\Omega}_{A_x/\C}^1,
\end{eqnarray*}
which shows that $\tau_{(z_1)}(\vp_{f^*}) = - v_1 \cdot dz_2 \neq 0$ 
and hence $C$ is of type I relative to $f$. 
\qed \par\medskip
Using this lemma we complete the proof of Theorem \ref{thm:main3}
\par\medskip\noindent
{\it Proof of Theorem $\ref{thm:main3}$}. 
Assume that the $2$-form $\omega$ has no pole of order $\nu_C(f)$ 
along $C$. 
In view of Lemma \ref{lem:typeII} the theorem is proved if 
$k > l$ is shown in (\ref{eqn:expand}). 
Around $x$ we express $\om$ as 
\begin{equation} \label{eqn:om}
\om = \a \cdot dz_1 \wedge dz_2 \qquad \mathrm{with} \quad 
\a = \sum_{n = s}^{\infty} \a_n(z_2) \, z_1^n
\end{equation}
where $s \in \Z$ and $\a_n(z_2) \in \C(\!(z_2)\!)$ with $\a_s(z_2)$ 
not identically zero. 
Assume the contrary that $k \le l$. 
Then we have $\nu_{C}(f) = k$ by Lemma \ref{lem:typeII} and 
hence $k \neq - s$, since $\omega$ has a pole of order 
$- s$ along $C$ (a pole of negative order is a zero). 
In order to consider the area-preserving property 
$f_x^* \omega = \omega$, we calculate  $f_x^*\om := 
\a(f_x^*(z_1), \, f_x^*(z_2)) \cdot 
d(f_x^*(z_1)) \wedge d(f_x^*(z_2))$. 
Considering the Laurent expansion of 
$\a(f_x^*(z_1), \, f_x^*(z_2))$ in $z_1$, we have 
\[
\begin{array}{rcl}
\a(f_x^*(z_1),f_x^*(z_2)) &=& 
\displaystyle \sum_{n=s}^{\infty} \a_n(z_2 + z_1^l \cdot f_2) 
\cdot (z_1+ z_1^k \cdot f_1)^n \\[6mm] 
&=& \displaystyle \sum_{i,j=0}^{\infty} \sum_{n = s}^{\infty} 
{n \choose i} \dfrac{\a_n^{(j)}(z_2)}{j!} 
\cdot f_1^i \cdot f_2^j \cdot z_1^{n+(k-1)i+lj} \\[6mm]
&=& \displaystyle \sum_{n=s}^{s+k-1} \a_n(z_2) \cdot z_1^n 
+ s \cdot \a_s (z_2) \cdot f_1 \cdot z_1^{s+k-1} + O(z_1^{s+k}), 
\end{array}
\]
where on the righthand side the first terms come from the indices 
$(i,j) = (0,0)$, $s \le n \le s+k-1$, the second term from 
$(i,j) = (1,0)$, $n = s$, and the $O(z_1^{s+k})$-term from the 
remaining indices, since we are assuming that $k-1 < l$. 
Similarly it follows from $1 \le k \le l$ that 
\[
\begin{array}{rcl}
d(f_x^*(z_1)) \wedge d(f_x^*(z_2)) 
&=& \{\, (1 + k \cdot z_1^{k-1} \cdot f_1 + z_1^k \cdot f_{1{z_1}})
\cdot (1 + z_1^l \cdot f_{2{z_2}}) \\[4mm]
& & \phantom{\{\,(1} - z_1^k \cdot f_{1{z_2} \,} \cdot
(l \cdot z_1^{l-1} \cdot f_2 + z_1^l \cdot f_{2{z_2}}) \} \cdot 
dz_1 \wedge dz_2 \\[4mm]
&=& \{1 + k z_1^{k-1} \cdot f_1 + O(z_1^k)\} \cdot 
dz_1 \wedge dz_2, 
\end{array}
\]
where $f_{i{z_j}}$ is the partial derivative of 
$f_i$ with respect to $z_j$. 
These calculations and (\ref{eqn:om}) yield 
\[
\begin{array}{rcl}
f_x^*\om  
&=& \left\{ \displaystyle \sum_{n=s}^{s+k-1} \a_n(z_2) \cdot z_1^n + 
(s+k) \cdot \a_s (z_2) \cdot f_1 \cdot z_1^{s+k-1} + O(z_1^{s+k}) 
\right\} \cdot dz_1 \wedge dz_2, \\[5mm]
\om 
&=& \left\{ \displaystyle \sum_{n=s}^{s+k-1} \a_n(z_2) \cdot z_1^n 
+ O(z_1^{s+k}) \right\} \cdot dz_1 \wedge dz_2. 
\end{array}
\]
Comparing the coefficients of $z_1^{s+k-1} \cdot dz_1 \wedge dz_2$ 
in the area-preserving condition $f_x^* \om = \om$, we have 
$(s+k) \cdot \a_s(z_2) \cdot f_1(0,z_2) = 0$ in $\C(\!(z_2)\!)$, 
but this contradicts the fact that $s + k \neq 0$ in $\Z$ and 
$\a_s(z_2) \neq 0$, $f_1(0,z_2) \neq 0$ in $\C(\!(z_2)\!)$. 
Therefore actually we have $k > l$ and hence $C$ is of type I\!I 
relative to $f$ by Lemma \ref{lem:typeII}. 
The proof is complete. \qed
\begin{remark} \label{rem:pre}
If $C$ is a fixed curve along which $\omega$ has a pole of order 
$\nu_C(f)$, then $C$ may be of type I. 
For example, consider the birational map $f : \P^2 \to \P^2$ mentioned 
in Remark \ref{rem:ind2}. 
Notice that $f$ preserves the meromorphic $2$-form 
$\omega := z_1^{-1} dz_1 \wedge dz_2$. 
The curve $C := \{z_1 = 0\}$ is a fixed curve of type I with 
index $\nu_C(f) = 1$, along which $\omega$ has a pole of order 
$\nu_C(f) = 1$. 
\end{remark}
\par 
In Section \ref{sec:intro} a discussion is made rather roughly 
as to the inapplicability of some generalized fixed formulas 
other than Saito's formula (\ref{eqn:formula}). 
We take this occasion to restate it more precisely. 
\begin{remark} \label{rem:index}
The generalized fixed point formulas cited in Section 
\ref{sec:intro} other than formula (\ref{eqn:formula}) are valid when 
$X$ is a complex manifold and $f : X \to X$ is a holomorphic map 
such that each connected component $Y$ of the fixed point set is a 
non-degenerate submanifold, that is, all the eigenvalues of 
the normal map $d^Nf : NY \to NY$ are different from $1$, where 
$d^Nf$ is the map induced from the tangent map $df : TX \to TX$ to 
the normal bundle $NY := T_{Y}X/TY$. 
Under this non-degeneracy condition, if $X$ is a surface, 
then Lemma \ref{lem:typeII} readily shows that any fixed curve of 
$f$ is of type I. 
On the other hand, if $f$ is an area-preserving surface map, 
then any fixed curve $C$ must be degenerate, because the 
normal map $d^Cf$ becomes identity. 
In fact we know from Theorem \ref{thm:main3} that $C$ is 
of type I\!I. 
Thus the area-preserving surface maps are completely outside 
the reach of the usual generalized fixed point formulas other 
than formula (\ref{eqn:formula}). 
\end{remark}
\section{Isolated Periodic Points} \label{sec:iso-per} 
In this section, based on the fundamental results in 
Theorems \ref{thm:ssaito}, \ref{thm:main}, \ref{thm:main2}, 
and \ref{thm:prime}, we are concerned with the isolated 
periodic points of a birational surface map $f$. 
In view of Theorem \ref{thm:main3} it is reasonable to 
assume that $f$ has no periodic curves of type I. 
Under this assumption we establish a certain periodic point 
formula (Theorem \ref{thm:formula}), a Shub-Sullivan type result 
(Lemma \ref{lem:SS2}) and a refined version of 
Theorem \ref{thm:main4} (Theorem \ref{thm:main4dash}). 
\par 
To this end we need to introduce some terminology and notation. 
Let $X_0^{c}(f)$ be the set of all non-isolated fixed points 
of $f$, namely, the set of all points that lie on some fixed 
curve of $f$, and let $X_0^{i}(f)$ be the set of all isolated 
fixed points of $f$, that is, the complement of $X_0^{c}(f)$ 
in $X_0(f)$. 
Then $\mathrm{Per}_n^{i}(f) := X_0^{i}(f^n)$ stands for 
the set of all isolated periodic points of $f$ with (not 
necessarily prime) period $n$ and its cardinality counted 
with multiplicity is defined by 
\begin{equation} \label{eqn:card-iso}
\# \mathrm{Per}_n^{i}(f) := \sum_{x \in \mathrm{Per}_n^{i}(f)} 
\nu_x(f^n). 
\end{equation}
Here a remark on notation: $\#$ is used to denote the 
cardinality counted with multiplicity or, in other words, 
the weighted cardinality, while $\mathrm{Card}$ is reserved for 
the cardinality without multiplicity taken into account. 
We denote by $P(f)$ the set of all positive integers that 
arise as the prime period of some irreducible periodic curve 
of $f$. 
For each $n \in \N$, let $P_n(f)$ be the set of all elements 
$k \in P(f)$ that divides $n$. 
Note that $P_n(f)$ is a finite set for every $n \in \N$, 
while $P(f)$ may or may not be finite. 
\begin{remark} \label{rem:P(f)} 
The map $f$ admits infinitely many irreducible periodic curves 
if and only if the set $P(f)$ is infinite, because the number 
of irreducible periodic curves of any given prime period is 
finite, provided that $f$ is nontrivial. 
(Recall that $f$ is always assumed to be nontrivial.) 
\end{remark} 
\par 
Moreover, for each $k \in P(f)$ we denote by $\mathrm{PC}_k(f)$  
the set of all irreducible periodic curves of $f$ with prime 
period $k$. 
There is then a direct sum decomposition: 
\begin{equation} \label{eqn:decomp1}
X_1(f^n) = \coprod_{k \in P_n(f)} \mathrm{PC}_k(f). 
\end{equation}
Given any $k \in P(f)$, let $C_k(f)$ be the (possibly 
reducible) curve in $X$ defined to be the union of all curves 
in $\mathrm{PC}_k(f)$. 
Then there exists a decomposition: 
\begin{equation} \label{eqn:decomp2}
X_0(f^n) = \mathrm{Per}_n^{i}(f) \amalg 
\bigcup_{k \in P_n(f)} C_k(f). 
\end{equation}
Finally, for each $k \in P(f)$, let $\xi_k(f)$ be the number 
defined by 
\begin{equation} \label{eqn:xik}
\xi_k(f) := \sum_{x \in C_k(f)} \nu_x(f^k) + 
\sum_{C \in \mathrm{PC}_k(f)} \tau_{C} \cdot \nu_{C}(f^k). 
\end{equation}
\par 
With these preliminaries, under the absence of periodic curves 
of type I, Saito's fixed point formula (\ref{eqn:formula}) 
is applied to the iterates $f^n$ to yield the following 
periodic point formula. 
\begin{theorem} \label{thm:formula} 
Let $f : X \to X$ be an AS birational map of a smooth 
projective surface $X$. 
If $f$ has no periodic curve of type I, then we have for any 
$n \in \N$, 
\begin{equation} \label{eqn:formula2} 
L(f^n) = \# \mathrm{Per}_n^{i}(f) + \sum_{k \in P_n(f)} \xi_k(f).  
\end{equation}
\end{theorem}
\proof
By assumption, the map $f^n$ has no fixed curves of type I, 
that is, $X_1(f^n) = X_{I\!I}(f^n)$ for any $n \in \N$. 
Then, by Theorem \ref{thm:prime}, if $k$, $l \in P(f)$ are 
distinct then $C_k(f)$ and $C_l(f)$ are disjoint.  
Thus (\ref{eqn:decomp2}) becomes the direct sum decomposition: 
\begin{equation} \label{eqn:decomp3}
X_0(f^n) = \mathrm{Per}_n^{i}(f) \amalg 
\coprod_{k \in P_n(f)} C_k(f). 
\end{equation}
Since $f$ is assumed to be AS, Remark \ref{rem:separate} 
implies that the fixed point formula (\ref{eqn:formula}) can 
be applied to all iterates $f^n$ $(n \in \N)$. 
In view of the direct sum decompositions 
(\ref{eqn:decomp1}) and (\ref{eqn:decomp3}) and the equality 
$X_1(f^n) = X_{I\!I}(f^n)$, the formula (\ref{eqn:formula}) 
is rewritten as 
\[
\begin{array}{rcl}
L(f^n) &=& \displaystyle \sum_{x \in X_0(f^n)} \nu_x(f^n) + 
\displaystyle \sum_{C \in X_1(f^n)} \tau_{C} \cdot \nu_{C}(f^n) 
\\[7mm]
&=& \displaystyle \sum_{x \in \mathrm{Per}_n^{i}(f)} \nu_x(f^n) + 
\displaystyle \sum_{k \in P_n(f)} \sum_{x \in C_k(f)} \nu_x(f^n) + 
\displaystyle \sum_{k \in P_n(f)} \sum_{C \in \mathrm{PC}_k(f)} 
\tau_{C} \cdot \nu_{C}(f^n) \\[7mm]
&=& \# \mathrm{Per}_n^{i}(f) + 
\displaystyle \sum_{k \in P_n(f)} 
\left\{ \sum_{x \in C_k(f)} \nu_x((f^k)^{n/k}) + 
\sum_{C \in \mathrm{PC}_k(f)} \tau_{C} \cdot \nu_{C}((f^k)^{n/k}) 
\right\}. 
\end{array}
\]
Here we note that $n/k \in \N$ for any $k \in P_n(f)$, any 
$x \in C_k(f)$ passes through a fixed curve of type I\!I of 
$f^k$ and any $C \in \mathrm{PC}_k(f)$ is a fixed curve of 
type I\!I of $f^k$. 
Thus Theorem \ref{thm:main} implies that $\nu_x((f^k)^{n/k}) 
= \nu_x(f^k)$ and  $\nu_{C}((f^k)^{n/k}) = \nu_{C}(f^k)$. 
Hence we have 
\[
\begin{array}{rcl}
L(f^n) &=& 
\# \mathrm{Per}_n^{i}(f) + 
\displaystyle \sum_{k \in P_n(f)} 
\left\{ \sum_{x \in C_k(f)} \nu_x(f^k) + 
\sum_{C \in \mathrm{PC}_k(f)} \tau_{C} \cdot \nu_{C}(f^k) 
\right\} \\[8mm]
&=& \# \mathrm{Per}_n^{i}(f) + 
\displaystyle \sum_{k \in P_n(f)} \xi_k(f), 
\end{array}
\]
where (\ref{eqn:xik}) is used in the last line. 
This proves the theorem. \qed \par\medskip 
We need a bit more terminology. 
Let $\mathrm{Per}^{i}(f)$ denote the set of all isolated 
periodic points of $f$, that is, the union of 
$\mathrm{Per}_n^{i}(f)$ over all $n \in \N$. 
We make the following definition. 
\begin{definition} \label{def:ab-iso} 
An isolated periodic point $x \in \mathrm{Per}^{i}(f)$ is said 
to be {\sl absolutely isolated} if $x$ is an isolated fixed point 
of $f^n$ for any period $n$ of $x$. 
Otherwise $x$ is said to be {\sl conditionally isolated}. 
Denote by $\mathrm{Per}^{ai}(f)$ and $\mathrm{Per}^{ci}(f)$ the 
set of all absolutely isolated periodic points and the set of all 
conditionally isolated periodic points of $f$ respectively. 
Then we have 
\[
\mathrm{Per}^{i}(f) = 
\mathrm{Per}^{ai}(f) \amalg \mathrm{Per}^{ci}(f). 
\]
For any $x \in \mathrm{Per}^{ci}(f)$ there exists a period 
$n \in \N$ of $x$ relative to $f$ such that $x$ is a non-isolated 
fixed point of $f^n$. 
The smallest such $n$ is called the {\sl secondary period} of 
$x$ relative to $f$. 
Put 
\[
\begin{array}{rclrcl}
\mathrm{Per}_n^{ai}(f) &:=& 
\mathrm{Per}^{ai}(f) \cap \mathrm{Per}_n^{i}(f), \qquad & 
\# \mathrm{Per}_n^{ai}(f) &:=& \displaystyle
\sum_{x \in \mathrm{Per}_n^{ai}(f)} \nu_x(f^n), \\[7mm]
\mathrm{Per}_n^{ci}(f) &:=& 
\mathrm{Per}^{ci}(f) \cap \mathrm{Per}_n^{i}(f), \qquad & 
\# \mathrm{Per}_n^{ci}(f) &:=& \displaystyle 
\sum_{x \in \mathrm{Per}_n^{ci}(f)} \nu_x(f^n).  
\end{array} 
\]
\end{definition}
\begin{remark} \label{rem:ab-iso} 
The secondary period $m$ of $x \in \mathrm{Per}^{ci}(f)$ is a 
strictly greater multiple of its prime period $n$. 
Indeed, if $m = n$ then $x$ would be a non-isolated fixed point
of $f^{ln}$ for any $l \in \N$, contradicting the assumption 
that $x$ is an isolated periodic point. 
Write $m = kn$ with $k \in \N_{\ge 2}$. 
Then $x$ is an isolated fixed point of $f^{ln}$ if and only if 
$l$ is not divisible by $k$. 
\end{remark} 
\par 
Shub and Sullivan \cite{SS} state their result just as in 
Theorem \ref{thm:SS}, but a careful check of their proof 
shows that their result is valid in the following more 
general form. 
\begin{theorem} \label{thm:SS2} 
Let $f : X \to X$ be a $C^1$-map of a smooth manifold $X$ and 
$x \in X$ an isolated fixed point of $f$. 
Let $N_x(f)$ be the set of all $n \in \N$ such that $x$ is an 
isolated fixed point of $f^n$. 
Then the indices $\nu_x(f^n)$ are bounded as a function of 
$n \in N_x(f)$. 
\end{theorem}
A combination of Theorems \ref{thm:main} and \ref{thm:SS2} 
leads to the following lemma. 
\begin{lemma} \label{lem:SS2} 
Assume that $f$ has no periodic curves of type I. 
Given any $x \in \mathrm{Per}^{i}(f)$, let 
$n$ be the prime period of $x$ relative to $f$. 
Then the indices $\nu_x(f^{ln})$ are bounded as a function 
of $l \in \N$. 
\end{lemma}
\proof 
First, if $x \in \mathrm{Per}^{ai}(f)$, then the lemma is proved 
by applying Theorem \ref{thm:SS} to the map $f^n$. 
Next consider the case $x \in \mathrm{Per}^{ci}(f)$. 
Let $m$ be the secondary period of $x$ relative to $f$ 
and put $m = k n$ with $k \in \N_{\ge 2}$. 
Note that $N_x(f^n) = \{\, l \in \N\,|\, \mbox{$l$ is not 
divisible by $k$}\,\}$ by Remark \ref{rem:ab-iso}. 
Then Theorem \ref{thm:SS2} applied to $f^n$ implies that 
$\nu_x(f^{ln})$ is bounded for $l \in N_x(f^n)$. 
For $l$ divisible by $k$, we can use Theorem \ref{thm:main} 
since $x$ is a fixed point of $f^m = f^{kn}$ through which 
a fixed curve of type I\!I of $f^m$ passes. 
The proof is complete. \qed \par\medskip
There is a relation between the number of periodic curves and 
that of conditionally isolated periodic points, as is stated 
in the following lemma. 
\begin{lemma} \label{lem:cond-iso} 
If $f$ has at most finitely many irreducible periodic 
curves or equivalently if the set $P(f)$ is finite 
$($see Remark $\ref{rem:P(f)}$$)$, then the set 
$\mathrm{Per}^{ci}(f)$ is also finite. 
\end{lemma}
\proof 
First we show that any conditionally isolated periodic point 
$x \in \mathrm{Per}^{ci}(f)$ is an intersection point of two or 
more distinct irreducible periodic curves of $f$. 
Let $n$ and $m$ be the prime period and the secondary period 
of $x$ relative to $f$, respectively. 
By Remark \ref{rem:ab-iso} there exists an integer $k \ge 2$ 
such that $m = kn$. 
Let $C$ be an irreducible fixed curve of $f^m$ passing 
through $x$. 
Since $x \in X_0(f^n)$, one has either $x \in X_0^{\circ}(f^n)$ 
or $x \in X_0^{\circ}(f^{-n})$ (see Definition \ref{def:fixed}). 
In the former case, for each $1 \le l \le k-1$, let 
$C_l := f^{ln}(C)$ be the strict transform of $C$ by $f^{ln}$. 
In the latter case we consider $C_l := f^{-ln}(C)$ instead. 
Then $C_l$ is an irreducible fixed curve of $f^m$ passing 
through $x$, but different from $C$. 
Thus $x \in C \cap C_l$ is an intersection point of the 
distinct irreducible periodic curves $C$ and $C_l$ of $f$. 
So the claim is verified. 
\par 
We proceed to the proof of the lemma. 
Assume that $f$ has at most finitely many irreducible 
periodic curves. 
Then the set of all intersection points of all pairs of 
distinct irreducible periodic curves is also finite. 
By what is proved in the last paragraph, 
$\mathrm{Per}^{ci}(f)$ is a subset of this finite set. 
Hence the set $\mathrm{Per}^{ci}(f)$ is also finite. 
\qed \par\medskip
Formula (\ref{eqn:formula2}) shows that the weighted cardinalities 
$\# \mathrm{Per}_n^{i}(f)$ are controlled by the magnitudes 
of the Lefschetz numbers $L(f^n)$ and the sets $P_n(f)$ of 
prime periods of periodic curves. 
In terms of the first dynamical degree $\lambda(f)$ the behavior 
of $L(f^n)$ is described as follows. 
\begin{lemma} \label{lem:exponential} 
If $f : X \to X$ is an AS bimeromorphic map on a compact K\"{a}hler 
surface $X$ with the first dynamical degree $\lambda(f) > 1$, then 
the Lefschetz numbers $L(f^n)$ admits the estimate: 
\[
|L(f^n) - \lambda(f)^n| \le \left\{ 
\begin{array}{ll}
O(1) \qquad &(\mbox{if $X \sim$ no complex $2$-torus}), \\[2mm] 
4 \, \lambda(f)^{n/2} + O(1) \quad &(\mbox{if $X \sim$ 
a complex $2$-torus}), 
\end{array} \right. 
\]
where $X \sim Y$ indicates that $X$ is bimeromorphically 
equivalent to $Y$. 
\end{lemma} 
\proof 
In what follows we use the notation: $h^i := \dim_{\C} H^i(X)$, 
$h^{i,j} := \dim_{\C} H^{i,j}(X)$ and 
$t_n^{i,j} := \mathrm{Tr}[\,(f^n)^* : H^{i,j}(X) \to H^{i,j}(X)\,]$. 
Note that $L(f^n) = \sum_{i,j} (-1)^{i+j} \, t_n^{i,j}$ and $t_n^{i,j}$ 
is the $n$-th power sum of the eigenvalues of $f^* : H^{i,j}(X) \carl$. 
First, for $i = 0, 2$, the induced map $(f^n)^* : H^{i,i}(X) \cong 
\C \carl$ is the identity since $f^n$ is a bimeromorphic map. 
Next we consider the case $(i,j) = (1, 1)$. 
Because $f$ is assumed to be AS, we have 
$(f^n)^* = (f^*)^n : H^{1,1}(X) \to H^{1,1}(X)$. 
It follows from \cite[Theorem 0.3]{DF} that $f^*|_{H^{1,1}(X)}$ has 
a simple eigenvalue $\lambda(f) > 1$ together with all the remaining 
eigenvalues in the closed unit disk in $\C$. 
This shows that 
\begin{equation} \label{eqn:t11n}
t^{1,1}_n = \lambda(f)^n + O(1). 
\end{equation}
\par 
It is well known that the first dynamical degree $\lambda(f)$ and 
the eigenvalues of $(f^n)^* : H^{i,j}(X) \carl$ with $(i,j) \neq 
(1,1)$ are invariant under meromorphic conjugation of $(X, f)$ 
(see \cite[Proposition 1.18]{DF} and \cite[page 34]{BHPV}). 
Hence for any $(i,j)$ the number $t_n^{i,j}$ is also invariant 
except for the $O(1)$-term of (\ref{eqn:t11n}) in the case 
$(i,j) = (1,1)$. 
So, in any case, every $t_n^{i,j}$ is invariant up to $O(1)$-term. 
Thus we may assume from the beginning that $(X, f)$ is of the 
canonical form (3) or (4) in Theorem \ref{thm:DF}. 
We make a case-by-case check. 
If $X$ is either a rational surface or an Enriques surface, 
then $h^1 = h^{2,0} = h^{0,2} = h^3 = 0$ and so 
(\ref{eqn:t11n}) implies that 
\begin{equation} \label{eqn:estimate1}
L(f^n) = t^{1,1}_n + 2 = \lambda(f)^n + O(1). 
\end{equation}
Next consider the case where $f$ is an automorphism of a K3 
surface $X$. 
Then we have $h^1 = h^3 = 0$ and 
$H^{2,0}(X) = \C\, \eta$ and $H^{0,2}(X) = \C \, \bar{\eta}$, 
where $\eta$ is a nowhere vanishing holomorphic $2$-form on $X$. 
There exists a constant $\delta \in \C^{\times}$ such that 
$f^*|_{H^{2,0}(X)}$ and $f^*|_{H^{0,2}(X)}$ are the scalar 
multiplications by $\delta$ and by $\bar{\delta}$ respectively. 
Since the automorphism $f$ preserves the volume form 
$\eta \wedge \bar{\eta}$, we have 
$\eta \wedge \bar{\eta} = f^* \eta \wedge \bar{\eta} = 
(\delta \eta) \wedge (\bar{\delta} \bar{\eta}) = 
|\delta|^2 \eta \wedge \bar{\eta}$ and $|\delta|=1$. 
Again (\ref{eqn:t11n}) yields: 
\begin{equation} \label{eqn:estimate2}
L(f^n) = t^{1,1}_n + \delta^n + \bar{\delta}^n +2 = 
\lambda(f)^n + O(1). 
\end{equation}
\par 
\begin{table}[t]
\begin{center}
\begin{tabular}{|c|c|}
\hline
\vspace{-3mm} &  \\
$(i,j)$ & the eigenvalues of $f^* : H^{i,j}(X) \carl$ \\[1mm]
\hline
\hline
\vspace{-3mm} &  \\
$(1, 0)$ & $\delta_1 = \delta$, \qquad 
$\delta_2 = \ve \delta^{-1}$ \\[1mm]
\hline
\vspace{-3mm} &  \\
$(0, 1)$ & $\bar{\delta}_1 = \bar{\delta}$, \qquad 
$\bar{\delta}_2 = \bar{\ve} \bar{\delta}^{-1}$  \\[1mm]
\hline
\vspace{-3mm} &  \\
$(2, 0)$ & $\delta_1 \delta_2 = \ve$  \\[1mm]
\hline
\vspace{-3mm} &  \\
$(0, 2)$ & $\bar{\delta}_1 \bar{\delta}_2 = \bar{\ve}$  \\[1mm]
\hline
\vspace{-3mm} &  \\
$(2, 1)$ & $\delta_1 \delta_2 \bar{\delta}_1 = \ve \bar{\delta}$, 
\qquad 
$\delta_1 \delta_2\bar{\delta}_2 = \bar{\delta}^{-1}$  \\[1mm]
\hline
\vspace{-3mm} &  \\
$(1, 2)$ & $\delta_1\bar{\delta}_1\bar{\delta}_2 = \bar{\ve}\delta$, 
\qquad 
$\delta_2 \bar{\delta}_1 \bar{\delta}_2 = \delta^{-1}$ \\[1mm]
\hline
\hline
\vspace{-3mm} &  \\
$(1, 1)$ & $|\delta_1|^2 = |\delta|^2$, \quad 
$\delta_1 \bar{\delta}_2 = \bar{\ve} (\delta/\bar{\delta})$, 
\quad 
$\delta_2 \bar{\delta}_1 = \ve (\bar{\delta}/\delta)$, \quad 
$|\delta_2|^2 = |\delta|^{-2}$  \\[1mm]
\hline 
\end{tabular}
\end{center}
\caption{The eigenvalues of $f^* : H^{i,j}(X) \carl$ for 
a complex $2$-torus $X$}
\label{tab:eigenvalue}
\end{table}
Finally we treat the case where $f$ is an automorphism 
of a complex $2$-torus $X = \C^2/\varGamma$ with 
$\varGamma \cong \Z^4$ being a lattice in $\C^2$. 
The map $f$ lifts to an affine automorohism 
\[
F : \C^2 \to \C^2, \quad 
(z_1, \, z_2) \mapsto (a_{11} z_1 + a_{12} z_2 + b_1, \, 
a_{21} z_1 + a_{22} z_2 + b_2) 
\]
through the canonical projection $\C \to \C/\varGamma$. 
The determinant of the matrix $A := (a_{ij}) \in M_2(\C)$ 
satisfies $|\det A|= 1$, since $F$ acts on $\varGamma$ 
bijectively. 
Let $\delta_1$ and $\delta_2$ be the eigenvalues of 
$A$ with $|\delta_1| \ge |\delta_2|$. 
Note that $|\det A| = 1$ implies 
$|\delta_1\delta_2| = 1$. 
If we put $\delta := \delta_1$ and $\ve := \delta_1 \delta_2$, 
then $|\delta| \ge 1$ and $|\ve| = 1$. 
The action $f^* : H^{1,0}(X) = \C \, dz_1 \oplus 
\C \, dz_2 \carl$ is represented by the matrix $A$ and so 
has the eigenvalues $\delta_1 = \delta$ and $\delta_2 = 
\ve \delta^{-1}$. 
In a similar manner, using the representations 
\[
H^{i,j}(X) = \bigoplus_{k_1 < \cdots < k_i} 
\bigoplus_{l_1 < \cdots < l_j} \C \, dz_{k_1}\wedge \cdots 
\wedge dz_{k_i} \wedge d\bar{z}_{l_1} \wedge \cdots 
\wedge d\bar{z}_{l_j}, 
\]
the eigenvalues of $f^* : H^{i,j}(X) \carl$ 
are given as in Table \ref{tab:eigenvalue}. 
So the Lefschetz numbers 
$L(f^n) = \sum_{i,j} (-1)^{i+j} \, t^{i,j}_n$ are 
calculated as 
\[
L(f^n) = |\delta|^{2n} + |\delta|^{-2n} + 2 
- 2 \, \mathrm{Re}\, \{\, (1+\bar{\ve}^n) \delta^n + 
(1+\ve^n) \delta^{-n} - \ve^n (1 + (\bar{\delta}/\delta)^n) 
\,\}. 
\]
Since $|\delta| \ge 1$, the spectral radius of 
$f^*|_{H^{1,1}(X)}$ is given by $|\delta|^2$ and we have 
$\lambda(f) = |\delta|^2 > 1$. 
Using this relation in the formula above, we can easily 
obtain the estimate: 
\begin{equation} \label{eqn:estimate3}
|L(f^n) - \lambda(f)^n | < 4 \lambda(f)^{n/2} + 11 
\qquad (n \in \N). 
\end{equation}
Now the lemma follows from the estimates (\ref{eqn:estimate1}), 
(\ref{eqn:estimate2}) and (\ref{eqn:estimate3}), where the 
constant $11$ in (\ref{eqn:estimate3}) should be replaced by 
$O(1)$ if $X$ is a proper modification of a $2$-torus. 
\qed \par\medskip 
With these preliminaries we establish the following theorem. 
\begin{theorem} \label{thm:main4dash} 
Let $X$ be a smooth projective surface and $f : X \to X$ an AS 
birational map without periodic curves of type I. 
If $\lambda(f) > 1$ then $f$ has at most finitely many irreducible 
periodic curves, at most finitely many conditionally isolated 
periodic points, and infinitely many absolutely isolated periodic 
points. 
Moreover we have $\# \mathrm{Per}_n^{ci}(f) = O(1)$ and 
\begin{equation} \label{eqn:isoperiodic2}
|\# \mathrm{Per}_n^{ai}(f) - \lambda(f)^n| \le 
\left\{ \begin{array}{ll}
O(1) \qquad & (\mbox{if $X \sim$ no Abelian surface}), \\[2mm] 
4 \, \lambda(f)^{n/2} + O(1) \quad & 
(\mbox{if $X \sim$ an Abelian surface}), 
\end{array} \right. 
\end{equation}
where $X \sim Y$ indicates that $X$ is birationally equivalent 
to $Y$. 
\end{theorem} 
\proof 
By the assumption that $f$ has no periodic curves of type I, 
any irreducible periodic curve of it is of type I\!I. 
Since $\lambda(f) > 1$, the assertion (1) of Theorem \ref{thm:main2} 
implies that $f$ has at most finitely many irreducible periodic 
curves, namely, that $\mathrm{Card} \, P(f) < \infty$.  
Since $P_n(f)$ is a subset of $P(f)$ for any $n \in \N$, 
the second term of the righthand side of formula 
(\ref{eqn:formula2}) is bounded as a function of $n \in \N$. 
So the Lefschetz numbers $L(f^n)$ and the weighted cardinalities 
$\# \mathrm{Per}_n^{i}(f)$ behave in the same manner 
modulo a bounded function of $n$: 
\begin{equation} \label{eqn:card-pern}
\# \mathrm{Per}_n^{i}(f) = L(f^n) + O(1). 
\end{equation}
For each $x \in \mathrm{Per}^{i}(f)$ let $n_x$ denote the 
prime period of $x$ relative to $f$. 
By Lemma \ref{lem:SS2} there exists a constant $M_x < \infty$ 
such that $\nu_x(f^{l\cdot n_x}) \le M_x$ for all $l \in \N$. 
Since $P(f)$ is finite, Lemma \ref{lem:cond-iso} 
implies that $\mathrm{Per}^{ci}(f)$ is also finite. 
As $n/n_x \in \N$ for every $x \in \mathrm{Per}_n^{i}(f)$, 
we have 
\[
\begin{array}{rclcl}
\# \mathrm{Per}_n^{ci}(f) 
&:=& \displaystyle \sum_{x \in \mathrm{Per}_n^{ci}(f)} \nu_x(f^n) 
&=& \displaystyle \sum_{x \in \mathrm{Per}_n^{ci}(f)} 
\nu_x(f^{(n/n_x)\cdot n_x}) 
\\[7mm]
&\le& \displaystyle \sum_{x \in \mathrm{Per}_n^{ci}(f)} M_x 
&\le& \displaystyle\sum_{x \in \mathrm{Per}^{ci}(f)} M_x < \infty, 
\end{array}
\]
which leads to $\# \mathrm{Per}_n^{ci}(f) = O(1)$. 
Then this together with (\ref{eqn:card-pern}) yields 
\begin{equation} \label{eqn:per-ai} 
\# \mathrm{Per}_n^{ai}(f) = L(f^n) + O(1). 
\end{equation}
We show that the set $\mathrm{Per}^{ai}(f)$ is infinite. 
Assume the contrary that it is finite. 
Then the same estimate as above with $\mathrm{Per}_n^{ci}(f)$ 
replaced by $\mathrm{Per}_n^{ai}(f)$ yields 
\[
\# \mathrm{Per}_n^{ai}(f) 
\le \displaystyle \sum_{x \in \mathrm{Per}^{ai}(f)} M_x < \infty. 
\]
This estimate and formula (\ref{eqn:per-ai}) imply that the 
Lefschetz numbers $L(f^n)$ are bounded, but this 
contradicts Lemma \ref{lem:exponential}. 
Finally, formula (\ref{eqn:isoperiodic2}) follows from Lemma 
\ref{lem:exponential} and formula (\ref{eqn:per-ai}). 
\qed\par\medskip\noindent
{\it Proof of Theorem $\ref{thm:main4}$}. 
By Theorem \ref{thm:main3} and condition $(*)$, the map $f$ has 
no periodic curves of type I. 
Then Theorem \ref{thm:main4} is an immediate consequence of 
Theorem \ref{thm:main4dash}. \qed\par\medskip
There is a counterpart of Theorem \ref{thm:main4dash} for the 
case $\lambda(f) = 1$. 
\begin{proposition} \label{prop:main4dash} 
Let $f : X \to X$ be an AS birational map without periodic 
curves of type I. 
If $\lambda(f) = 1$ and the Lefschetz numbers $L(f^n)$ are 
unbounded, then either 
\begin{enumerate}
\item $f$ has at most finitely many irreducible periodic 
curves, at most finitely many conditionally isolated periodic 
points, and infinitely many absolutely isolated periodic points; 
or 
\item $f$ has infinitely many irreducible periodic curves and 
preserves a unique rational or elliptic fibration such that 
any irreducible periodic curve is contained in a fiber of 
the fibration. 
\end{enumerate}
\end{proposition}
\proof 
Since the Lefschetz numbers $L(f^n)$ are unbounded, 
$f^n$ is not isotopic to the identity for any $n \in \N$. 
Thus we are in case (1) or (2) of Theorem \ref{thm:DF}, so 
that $f$ preserves a unique rational or elliptic fibration. 
The remaining proof is similar to the proof of 
Theorem \ref{thm:main4dash}, again making use of 
Theorem \ref{thm:formula}, Lemmas \ref{lem:SS2} and 
\ref{lem:cond-iso}. 
The difference is to apply assertion (2) of 
Theorem \ref{thm:main2} instead of assertion (1) of the 
same theorem, and to apply Theorem \ref{thm:zero-deg} instead 
of Theorems \ref{thm:pos-kodaira} and \ref{thm:rational}. 
Details may be omitted. \qed
\section{An Example} \label{sec:example}
In order to illustrate our main theorems, we give an example 
of an AS birational map preserving a meromorphic $2$-form 
on a smooth projective rational surface. 
This example arises as a special case of a $4$-parameter family 
of dynamical systems on cubic surfaces derived from the nonlinear 
monodromy of the sixth Painlev\'e equation via the Riemann-Hilbert 
correspondence \cite{IIS,IU1}. \par 
Let $\overline{S}$ be the projective cubic surface in $\P^3$ 
defined by the homogeneous cubic equation: 
\[
Z_1Z_2Z_3  + Z_0(Z_1^2 +Z_2^2 +Z_3^2) - 
8 Z_0^2(Z_1 + Z_2 + Z_3) + 28 Z_0^3 = 0, 
\]
in homogeneous coordinates $Z = [Z_0:Z_1:Z_2:Z_3]$. 
It has a unique singularity at 
\[
q = [1:2:2:2], 
\]
which turns out to be a simple singularity of type $D_4$. 
The intersection of $\overline{S}$ with the plane $\{Z_0 = 0\}$ 
at infinity yields tritangent lines $L_i = \{ Z_0 = Z_i = 0 \}$ 
at infinity $(i = 1,2,3)$. 
Then the affine cubic surface $S := \overline{S} \sm L$ with 
$L:= L_1 \cup L_2 \cup L_3$ is given by the affine cubic 
equation: 
\[
g(z) := z_1z_2z_3  + z_1^2 +z_2^2 +z_3^2 - 
8 (z_1 + z_2 + z_3) + 28 = 0, 
\]
where $z_i := Z_i/Z_0$. 
Since this equation is quadratic in each variable $z_i$, the 
line through a point $z \in S$ parallel to the $z_i$-axis 
passes through a unique second point $\si_i(z) \in S$. 
Hence we have three involutive automorphisms $\si_i : S \to S$, 
which are written explicitly as 
\[
\si_i : (z_i, \, z_j, \, z_k) \mapsto 
(8 - z_i- z_j z_k, \, z_j, \, z_k) \qquad (i = 1,2,3),
\]
with $\{i,j,k\} = \{1,2,3\}$. 
Note that the singular point $q$ is a fixed point of the 
involutions $\si_i$. \par 
A natural (complex) area-form on $S$ is given by its 
Poincar\'{e} residue: 
\[
\omega_{S} := 
\dfrac{dz_1 \wedge dz_2 \wedge dz_3}{dg} \qquad 
\mbox{restricted on} \quad S \sm \{q \}. 
\]
The map $\si_i$ sends $\omega_{S}$ to its negative: 
$\si_i^* \omega_{S} = - \omega_{S}$. 
Moreover $\si_i$ extends to a birational map 
$\overline{\si}_i : \overline{S} \to \overline{S}$ and 
$\omega_{S}$ extends to a $2$-form $\omega_{\overline{S}}$ 
which is holomorphic on $S \sm \{q \}$ and meromorphic on 
$\overline{S}\sm \{ q \}$ with simple poles along the 
tritangent lines at infinity: 
$(\omega_{\overline{S}})_{\infty} = L_1 + L_2 + L_3$. 
We consider an $\omega_{\overline{S}}$-preserving birational 
map defined by 
\[
\overline{\si} := (\overline{\si}_1 \circ 
\overline{\si}_2 \circ \overline{\si}_3)^2 : 
\overline{S} \to \overline{S}. 
\]
For the same reason as in \cite[Lemma 16]{IU1}, $\overline{\si}$ 
has no periodic curves of any prime period. \par
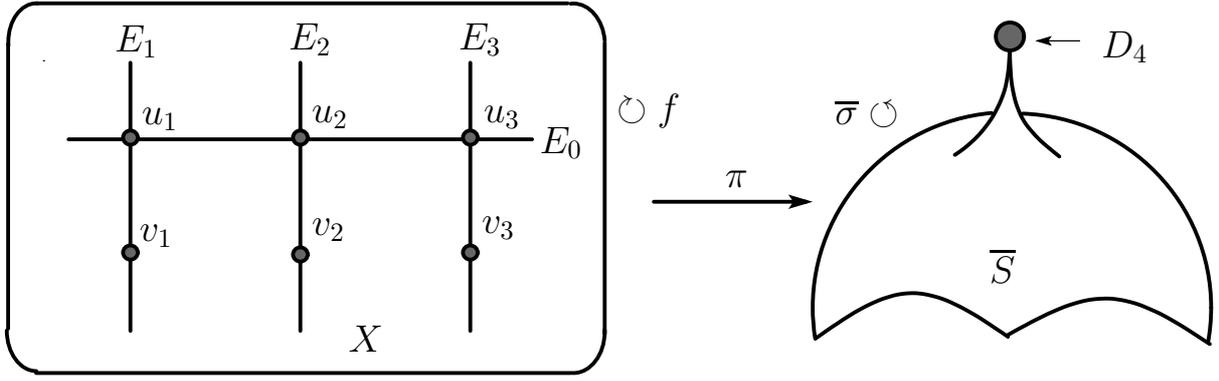
\begin{figure}[t]
\begin{center}
\unitlength 0.1in
\begin{picture}(62.36,19.30)(2.00,-20.20)
%
\special{pn 20}%
\special{pa 352 90}%
\special{pa 3104 90}%
\special{fp}%
%
\special{pn 20}%
\special{pa 352 2020}%
\special{pa 3120 2020}%
\special{fp}%
%
\special{pn 20}%
\special{pa 207 302}%
\special{pa 207 1788}%
\special{fp}%
%
\special{pn 20}%
\special{pa 3296 296}%
\special{pa 3296 1815}%
\special{fp}%
%
\special{pn 8}%
\special{ar 390 390 0 0  3.1415927 3.1967852}%
%
\special{pn 20}%
\special{ar 353 291 145 201  3.0438657 4.7606274}%
%
\special{pn 20}%
\special{ar 353 1798 153 212  1.5707963 3.1415927}%
%
\special{pn 20}%
\special{ar 3134 1805 162 215  6.2276868 6.2831853}%
\special{ar 3134 1805 162 215  0.0000000 1.6154593}%
%
\special{pn 20}%
\special{ar 3134 307 163 216  4.6227423 6.1786297}%
%
\special{pn 20}%
\special{pa 5394 263}%
\special{pa 5394 295}%
\special{pa 5393 328}%
\special{pa 5392 360}%
\special{pa 5390 392}%
\special{pa 5388 424}%
\special{pa 5384 456}%
\special{pa 5379 487}%
\special{pa 5373 519}%
\special{pa 5365 549}%
\special{pa 5354 580}%
\special{pa 5342 610}%
\special{pa 5328 639}%
\special{pa 5313 668}%
\special{pa 5296 696}%
\special{pa 5277 723}%
\special{pa 5257 749}%
\special{pa 5236 774}%
\special{pa 5213 798}%
\special{pa 5189 820}%
\special{pa 5165 841}%
\special{pa 5139 861}%
\special{pa 5113 879}%
\special{pa 5109 881}%
\special{sp}%
%
\special{pn 20}%
\special{pa 5394 263}%
\special{pa 5394 297}%
\special{pa 5394 332}%
\special{pa 5394 366}%
\special{pa 5395 400}%
\special{pa 5397 433}%
\special{pa 5399 466}%
\special{pa 5402 499}%
\special{pa 5407 531}%
\special{pa 5413 562}%
\special{pa 5420 593}%
\special{pa 5429 623}%
\special{pa 5440 651}%
\special{pa 5452 679}%
\special{pa 5467 706}%
\special{pa 5484 731}%
\special{pa 5504 755}%
\special{pa 5525 779}%
\special{pa 5547 802}%
\special{pa 5571 824}%
\special{pa 5596 845}%
\special{pa 5621 866}%
\special{pa 5647 887}%
\special{pa 5650 890}%
\special{sp}%
%
\special{pn 20}%
\special{ar 5394 1679 1021 1021  2.9830536 4.6192473}%
%
\special{pn 20}%
\special{ar 5422 1679 1014 1014  4.7504467 6.2831853}%
\special{ar 5422 1679 1014 1014  0.0000000 0.1784423}%
%
\special{pn 20}%
\special{pa 4387 1840}%
\special{pa 4416 1820}%
\special{pa 4445 1799}%
\special{pa 4474 1779}%
\special{pa 4503 1760}%
\special{pa 4532 1740}%
\special{pa 4561 1722}%
\special{pa 4589 1704}%
\special{pa 4618 1688}%
\special{pa 4647 1672}%
\special{pa 4676 1658}%
\special{pa 4705 1645}%
\special{pa 4734 1633}%
\special{pa 4763 1624}%
\special{pa 4792 1615}%
\special{pa 4821 1609}%
\special{pa 4850 1605}%
\special{pa 4879 1603}%
\special{pa 4908 1603}%
\special{pa 4937 1606}%
\special{pa 4966 1611}%
\special{pa 4995 1618}%
\special{pa 5024 1626}%
\special{pa 5053 1636}%
\special{pa 5082 1648}%
\special{pa 5111 1662}%
\special{pa 5140 1676}%
\special{pa 5169 1692}%
\special{pa 5199 1709}%
\special{pa 5228 1726}%
\special{pa 5257 1745}%
\special{pa 5286 1764}%
\special{pa 5315 1784}%
\special{pa 5344 1803}%
\special{pa 5373 1823}%
\special{pa 5384 1831}%
\special{sp}%
%
\special{pn 20}%
\special{pa 5384 1821}%
\special{pa 5414 1804}%
\special{pa 5444 1787}%
\special{pa 5475 1771}%
\special{pa 5505 1754}%
\special{pa 5535 1739}%
\special{pa 5565 1724}%
\special{pa 5595 1709}%
\special{pa 5625 1696}%
\special{pa 5655 1683}%
\special{pa 5685 1671}%
\special{pa 5715 1661}%
\special{pa 5745 1652}%
\special{pa 5775 1644}%
\special{pa 5805 1638}%
\special{pa 5835 1634}%
\special{pa 5865 1632}%
\special{pa 5894 1631}%
\special{pa 5924 1632}%
\special{pa 5953 1636}%
\special{pa 5983 1641}%
\special{pa 6012 1648}%
\special{pa 6042 1656}%
\special{pa 6071 1666}%
\special{pa 6100 1678}%
\special{pa 6129 1691}%
\special{pa 6159 1704}%
\special{pa 6188 1719}%
\special{pa 6217 1735}%
\special{pa 6246 1752}%
\special{pa 6275 1769}%
\special{pa 6304 1787}%
\special{pa 6333 1806}%
\special{pa 6362 1825}%
\special{pa 6391 1844}%
\special{pa 6420 1863}%
\special{pa 6429 1869}%
\special{sp}%
\put(39.2000,-10.5000){\makebox(0,0)[lb]{\large $\pi$}}%
%
\special{pn 20}%
\special{ar 5394 263 0 0  1.4141944 6.2831853}%
\special{ar 5394 263 0 0  0.0000000 1.2490458}%
%
\special{pn 20}%
\special{sh 0.600}%
\special{ar 5394 263 73 73  0.0000000 6.2831853}%
\put(52.9000,-15.5000){\makebox(0,0)[lb]{\large $\overline{S}$}}%
%
\special{pn 20}%
\special{pa 3550 1126}%
\special{pa 4338 1126}%
\special{fp}%
\special{sh 1}%
\special{pa 4338 1126}%
\special{pa 4271 1106}%
\special{pa 4285 1126}%
\special{pa 4271 1146}%
\special{pa 4338 1126}%
\special{fp}%
\put(19.7000,-19.1000){\makebox(0,0)[lb]{\large $X$}}%
\put(58.6900,-4.0000){\makebox(0,0)[lb]{\large $D_4$}}%
%
\special{pn 8}%
\special{pa 5755 286}%
\special{pa 5536 286}%
\special{fp}%
\special{sh 1}%
\special{pa 5536 286}%
\special{pa 5603 306}%
\special{pa 5589 286}%
\special{pa 5603 266}%
\special{pa 5536 286}%
\special{fp}%
%
\special{pn 20}%
\special{pa 520 800}%
\special{pa 2920 800}%
\special{fp}%
%
\special{pn 20}%
\special{pa 840 400}%
\special{pa 840 1800}%
\special{fp}%
\special{pa 2600 400}%
\special{pa 2600 1800}%
\special{fp}%
\special{pa 1720 400}%
\special{pa 1720 1800}%
\special{fp}%
\put(29.6000,-9.0000){\makebox(0,0)[lb]{\large $E_0$}}%
\put(7.6000,-3.7000){\makebox(0,0)[lb]{\large $E_1$}}%
\put(16.6000,-3.6000){\makebox(0,0)[lb]{\large $E_2$}}%
\put(25.4000,-3.6000){\makebox(0,0)[lb]{\large $E_3$}}%
%
\special{pn 20}%
\special{sh 0.600}%
\special{ar 840 791 38 37  0.0000000 6.2831853}%
%
\special{pn 20}%
\special{sh 0.600}%
\special{ar 1720 791 38 37  0.0000000 6.2831853}%
%
\special{pn 20}%
\special{sh 0.600}%
\special{ar 2600 791 38 37  0.0000000 6.2831853}%
%
\special{pn 20}%
\special{sh 0.600}%
\special{ar 2600 1391 38 37  0.0000000 6.2831853}%
%
\special{pn 20}%
\special{sh 0.600}%
\special{ar 1720 1401 38 37  0.0000000 6.2831853}%
%
\special{pn 20}%
\special{sh 0.600}%
\special{ar 840 1391 38 37  0.0000000 6.2831853}%
\put(9.0000,-7.5000){\makebox(0,0)[lb]{\large $u_1$}}%
\put(17.8000,-7.4000){\makebox(0,0)[lb]{\large $u_2$}}%
\put(26.7000,-7.4000){\makebox(0,0)[lb]{\large $u_3$}}%
\put(8.8800,-13.6000){\makebox(0,0)[lb]{\large $v_1$}}%
\put(17.8000,-13.3000){\makebox(0,0)[lb]{\large $v_2$}}%
\put(26.6000,-13.2000){\makebox(0,0)[lb]{\large $v_3$}}%
\put(33.6000,-7.4000){\makebox(0,0)[lb]{\large $\circlearrowright f$}}%
\put(44.9000,-7.4000){\makebox(0,0)[lb]{\large $\overline{\sigma} \circlearrowleft$}}%
\end{picture}%
\end{center}
\caption{Minimal resolution of singularities of type $D_4$} 
\label{fig:exceptional} 
\end{figure}

 
Let $\pi : (X, E) \to (\overline{S}, q)$ be a minimal 
desingularization of $\overline{S}$. 
Then its exceptional set $E$ consists of four irreducible 
components $E_0$, $E_1$, $E_2$, $E_3$ as depicted in 
Figure \ref{fig:exceptional}. 
We denote the strict transform of $L_i$ by the same symbol $L_i$. 
Then the pull-back $\omega_X := \pi^* \omega_{\overline{S}}$ 
turns out to be a meromorphic $2$-form on $X$ with simple poles 
along the tritangent lines at infinity: 
\[
(\omega_X)_{\infty} = L := L_1 + L_2 + L_3.  
\]
Note that $\omega_X$ is holomorphic and nondegenerate on 
$X \setminus L$ even around $E$. 
Let $p_i$ be the intersection point of $L_j$ and $L_k$ for 
$\{i,j,k\} = \{1,2,3\}$. 
It is easy see that the map $\overline{\sigma} : \overline{S} 
\to \overline{S}$ lifts to a birational map $f : X \to X$ such 
that $f^{-n}I(f) = \{p_3\}$ and $f^n I(f^{-1}) = \{p_1\}$ for 
any $n \in \N$ (see \cite[formula (52)]{IU1}). 
Hence $f$ is AS. 
We observe that $E_0$, $E_1$, $E_2$ and $E_3$ are irreducible 
fixed curves of $f$. 
There is no other periodic curves of $f$, since 
$\overline{\sigma}$ has no periodic curves as mentioned earlier. 
Thus under the notation of Section \ref{sec:APD} we have 
\begin{equation} \label{eqn:PC1}
P(f) = \{1\}, \qquad 
\mathrm{PC}_1(f) = \{E_0, E_1, E_2, E_3 \}, 
\qquad 
C_1(f) = E : = E_0 \cup E_1 \cup E_2 \cup E_3. 
\end{equation}
Note that $f$ preserves the meromorphic $2$-form $\om_X$. 
Since $\om_X$ is holomorphic in a neighborhood of $E$, 
Theorem \ref{thm:main3} implies that each fixed curve 
$E_i$ is of type I\!I relative to $f$. 
We wish to calculate the number $\xi_1(f)$ defined in 
(\ref{eqn:xik}). \par
Let $u_i$ denote the intersection point of $E_0$ and $E_i$ 
$(i = 1,2,3)$. 
We calculate the indices $\nu_{u_i}(f)$ and $\nu_{E_j}(f)$ 
$(j = 0,1,2,3)$. 
The induced endomorphism $f_{u_i}^* : A_{u_i} \to A_{u_i}$ 
is represented as 
\[
f_{u_i}^*(z_1) = z_1 + z_1^3 z_2 \, h_1(z), \quad 
f_{u_i}^*(z_2) = z_2 + z_1^2 z_2^2 \, h_2(z), 
\]
with some units $h_1$, $h_2 \in A_{u_i}^{\times}$, where 
$z_1$ and $z_2$ are local coordinates around $u_i$ such that 
$E_0 = \{z_1 = 0\}$ and $E_i = \{z_2 = 0\}$. 
By the definitions of $\ag(\sigma)$ and $\bg(\sigma)$ in 
Section \ref{sec:formula} and (\ref{eqn:omega}), 
\[
\begin{array}{rcl}
\ag(f_{u_i}^*) &=& (z_1^2z_2), \\[2mm] 
\bg(f_{u_i}^*) &=& (z_1 \, h_1(z), \, z_2 \, h_2(z)) = 
(z_1, \, z_2), \\[2mm]
\vp_{f_{u_i}^*} &=& z_2 h_2(z) \cdot dz_1 - z_1 h_1(z) \cdot dz_2 
\in \hat{\Om}_{A_{u_i}/\C}^1, 
\end{array}
\]
from which (\ref{eqn:delta}), (\ref{eqn:nup}), (\ref{eqn:mup}) and 
(\ref{eqn:nuC}) yield 
\[
\begin{array}{rcl}
\delta(f_{u_i}^*) &=& \dim_{\C} \, A_{u_1}/(z_1, \, z_2) = 1, \\[2mm]
\nu_{(z_1)}(f_{u_i}^*) &=& \nu_{E_0}(f) = 2, \qquad \\[2mm]
\nu_{(z_2)}(f_{u_i}^*) &=& \nu_{E_i}(f)=1, \\[2mm]
\mu_{(z_1)}(f_{u_i}^*) &=& \mu_{(z_2)}(f_{u_i}^*) = 1. 
\end{array}
\] 
By substituting these results into (\ref{eqn:nuA}) 
and (\ref{eqn:nux}), the index of $f$ at $u_i$ is given by 
\[
\nu_{u_i}(f) = \delta(f_{u_i}^*) + 
\sum_{k=1}^2 \nu_{(z_k)}(f_{u_i}^*) \cdot 
\mu_{(z_k)}(f_{u_i}^*) = 4. 
\]
\par 
An extra work shows that besides $u_1$, $u_2$, $u_3$, there are 
exactly three other points $x \in E$ such that $\nu_x(f) > 0$. 
More precisely, for each $i = 1,2,3$, there is a unique such 
point $v_i \in E_i \sm \{ u_i \}$, at which one has 
$\nu_{v_i}(f) = 2$ (see Figure \ref{fig:exceptional}). 
Summarizing these calculations, we have 
\begin{eqnarray} 
\nu_{x}(f) &=& \left\{
\begin{array}{rl}
4 \qquad &(x = u_i, \,\, i = 1,2,3), \\[2mm]
2 \qquad &(x = v_i, \,\, i = 1,2,3), \\[2mm]
0 \qquad &(\mbox{$x$ : any other point on $E$}), 
\end{array}
\right. \label{eqn:ex-index1} \\
\nonumber \\ 
\nu_{E_i}(f) &=& \left\{
\begin{array}{rl}
2 \qquad &(i=0), \\[2mm]
1 \qquad &(i=1,2,3). 
\end{array}
\right. \label{eqn:ex-index2} 
\end{eqnarray}
\par  
Since $\tau_{E_i} = - 2$ for each $i \in \{0,1,2,3\}$, 
substituting (\ref{eqn:PC1}), (\ref{eqn:ex-index1}) and 
(\ref{eqn:ex-index2}) into (\ref{eqn:xik}) yields 
\[
\xi_1(f) := \displaystyle \sum_{x \in C_1(f)} \nu_x(f) + 
\sum_{C \in \mathrm{PC}_1(f)} \tau_{C} \cdot \nu_C(f) 
= \displaystyle \sum_{x \in E} \nu_x(f) + 
\sum_{i = 0}^3 \tau_{E_i} \cdot \nu_{E_i}(f) = 8. 
\]
For every $n \in \N$ one has $P_n(f) = P(f) = \{1\}$ and thus 
the fixed point formula (\ref{eqn:formula2}) yields 
\[
L(f^n) = \# \mathrm{Per}_n^{i}(f) + \xi_1(f)
=  \# \mathrm{Per}_n^{i}(f) + 8. 
\]
On the other hand, we are able to show that 
$\lambda(f) = 9 + 4 \sqrt{5}$ and $L(f^n) = 
\lambda(f)^n + \lambda(f)^{-n} + 6$. 
Therefore we arrive at the following explicit formula for 
the number of isolated periodic points: 
\[
\# \mathrm{Per}_n^{i}(f) = 
(9 + 4 \sqrt{5})^n + (9 + 4 \sqrt{5})^{-n} - 2. 
\]
For this example all the isolated periodic points are 
absolutely isolated periodic points. 

\end{document}